\documentstyle{amlts}
\begin{document}
\annalsline{153}{2001}
\received{April 9, 1999}
\startingpage{661}
\def\bye{

 \end{document}}
 \font\tenrm=cmr10
\input amssym.def
\input amssym.tex
\def\ritem#1{\item[{\rm #1}]}

\catcode`\@=11
\font\twelvemsb=msbm10 scaled 1100
\font\tenmsb=msbm10
\font\ninemsb=msbm10 scaled 800
\newfam\msbfam
\textfont\msbfam=\twelvemsb  \scriptfont\msbfam=\ninemsb
  \scriptscriptfont\msbfam=\ninemsb
\def\msb@{\hexnumber@\msbfam}
\def\Bbb{\relax\ifmmode\let\next\Bbb@\else
 \def\next{\errmessage{Use \string\Bbb\space only in math
mode}}\fi\next}
\def\Bbb@#1{{\Bbb@@{#1}}}
\def\Bbb@@#1{\fam\msbfam#1}
\catcode`\@=12

 \catcode`\@=11
\font\twelveeuf=eufm10 scaled 1100
\font\teneuf=eufm10
\font\nineeuf=eufm7 scaled 1100
\newfam\euffam
\textfont\euffam=\twelveeuf  \scriptfont\euffam=\teneuf
  \scriptscriptfont\euffam=\nineeuf
\def\euf@{\hexnumber@\euffam}
\def\frak{\relax\ifmmode\let\next\frak@\else
 \def\next{\errmessage{Use \string\frak\space only in math
mode}}\fi\next}
\def\frak@#1{{\frak@@{#1}}}
\def\frak@@#1{\fam\euffam#1}
\catcode`\@=12

\newcommand{\bref}[1]{(\ref{#1})}
\newcommand{\fr}{{\frac{1}{2}}}
\renewcommand{\l}{{\lambda}}
\newcommand{\e}{{\varepsilon}}
\newcommand{\g}{{\gamma}}
\newcommand{\ld}{{\mbox{{\rm log}}\frac{1}{\delta}}}
\newcommand{\M}{{{\cal M}}}
\newcommand{\G}{{{\cal G}}}
\newcommand{\F}{{{\cal F}}}
\newcommand{\A}{{{\cal A}}}
\renewcommand{\r}{{\rho}}
\renewcommand{\d}{{\delta}}
\newcommand{\D}{{\Delta}}
\newcommand{\w}{{^{(1)}}}
\renewcommand{\t}{{^{(2)}}}
\newcommand{\p}{{\Phi\cdot\tilde\Phi}}
\renewcommand{\j}{{j^{\ast}}}
\newcommand{\s}{{\sigma}}
\newcommand{\su}{\mbox{supp}}
\renewcommand{\e}{{\varepsilon}}
\newcommand{\R}{{{\Bbb R}}}
\newcommand{\Z}{{{\Bbb Z}}}
\newcommand{\md}{{M_{\delta}}}
\newcommand{\ov}[1]{{\overline{#1}}}
\newcommand{\C}{{{\cal C}}}
\renewcommand{\O}{{{\cal O}}}
\newcommand{\J}{{{\cal J}}}
\newcommand{\W}{{\cal W}}
\newcommand{\Q}{{{\cal Q}}}
\newcommand{\B}{{{\cal B}}}
\newcommand{\lo}[1]{{\left(\log\frac{1}{\d}\right)^{#1}}}
\newcommand{\ts}{{\tilde{s}}}
\newcommand{\di}{\mbox{dist}}
\newcommand{\hf}{\mbox{\hfill$\square}}
\renewcommand{\c}{{\C}}
\newcommand{\0}{{{\cal O}}}
\newcommand{\gi}{{\gamma_i(s_i)}}
\newcommand{\x}[1]{{\gamma_{#1}(s_{#1})}}
\newcommand{\dgi}{{\dot{\gamma}_i(s_i)}}
\newcommand{\dgio}{{\dot{\gamma}_i(0)}}
\newcommand{\dg}{{\dot{\gamma}(s_i)}}
\newcommand{\dgo}{{\dot{\gamma}(0)}}
\newcommand{\gio}{{\gamma_i(0)}}
\renewcommand{\a}{{\A}}
\newcommand{\et}{{\sqrt{\frac{\e}{t}}}}
\setlength{\unitlength}{1cm}
\linethickness{.4mm}

\title{A sharp bilinear cone restriction estimate} 
\shorttitle{Sharp bilinear cone restriction estimate} 

 \author{Thomas Wolff}
\vglue-32pt
\centerline{\phantom{By Thomas Wolff}\quad$^{^\dagger}$}

 \institutions{California Institute of Technology, Pasadena, CA}

 \vglue12pt
The purpose of this paper is to prove an essentially sharp $L^2$ Fourier restriction
 estimate for light cones, of the type which is called bilinear in the recent literature. 

 Fix $d\geq 3$, denote variables in $\R^d$ by $(\ov{x},x_d)$ with $\ov{x}\in\R^{d-1}$,
 and let $\Gamma=\{ x:x_d=|\ov{x}| \mbox{ and } 1\leq x_d\leq 2\}$. Let $\Gamma_1$
 and $\Gamma_2$ be disjoint conical subsets, i.e.\ $$\Gamma_i=\{x\in\Gamma:\frac{\ov{x}}{x_d}\in\Omega_i\}$$
where $\Omega_i$ are disjoint closed subsets of the sphere $S^{d-2}$. Let $f$ and $g$ 
be two functions on $\Gamma$ whose supports are contained in  $\Gamma_1$ and
 $\Gamma_2$ respectively. We will prove the following estimate, where  $\s$ 
is surface measure on $\Gamma$, and
$\widehat{fd\s}$ is the  $\R^d$ Fourier transform:

\nonumproclaim{Theorem 1} If $p>1+\frac{2}{d}$ then 
\begin{equation}\|\widehat{fd\s}\widehat{gd\s}\|_p\leq C_{p, \Gamma_1,\Gamma_2}\|f\|_2\|g\|_2 . \label{5/3}
\end{equation}
\endproclaim
 
\vglue-12pt 

Bilinear estimates of this general type have been used by several authors; see in particular \cite{KM}. 
The
 estimate \bref{5/3} was 
formulated by Bourgain in \cite{B}, and it was proved in \cite{B}
 when 
$d=3$ and $p>2-\e$ for some $\e>0$,  the case $p=2$ being easier and 
implicit in \cite{Ba}.  Tao and Vargas \cite{TV}  recently obtained the explicit range
 $p>2-\frac{8}{121}$ when $d=3$, and noted that one can also obtain a range 
$p>2-\e_0$ in the four dimensional case. The range of $p$ in Theorem 1
is known to be best possible when $d=3$ except for the question of the endpoint --
see \cite{TV}, where   the conjecture that \bref{5/3} should hold for $d=3$ and  $p\geq\frac{5}{3}$ is
 attributed to  Machedon and Klainerman  -- and is similarly 
best possible in higher dimensions; see \cite{FK}.  

Although Theorem 1 is sharp of its type in any dimension, it is
 more satisfactory in low dimensions, since when $d$ is large the
 $L^2$ norms on the right hand side of \bref{5/3} are quite weak in
 comparison with other relevant norms and the exponent $1+\frac{2}{d}$ 
is only a small improvement on the exponent $1+\frac{2}{d-2}$ which
\vfill

\footnoterule
{\footnotesize $^\dagger$Thomas Wolff died tragically on July 31, 2000, after submitting this paper.}
\eject
\noindent  follows 
from the Strichartz inequality.  When $d=4$, Theorem 1 implies  (via a rescaling argument as in \cite{TVV}) a
statement analogous to a result of Barcelo \cite{Ba} for the three dimensional
 case:

\nonumproclaim{{C}orollary} When $d=4$ the restriction of the Fourier transform to $\Gamma$
 defines a
bounded operator from $L^p$ to $L^p(\Gamma)$ for any $p<\frac{3}{2}${\rm .}
\endproclaim

The range of $p$ here is again sharp. It should be pointed out that  the geometric
 information needed for our results is  simpler than what is likely to
 be needed either to solve the restriction problem for $S^2$, or to solve some of
 the other outstanding problems concerning the cone such as the 
multiplier problem and
 local smoothing, even in the $2+1$-dimensional case. 
 On the 
other hand, there are very few hypersurfaces for which a  sharp restriction theorem is
 known, and  the approach below may  be useful in connection with
 the sphere as well, insofar as it is possible to consider the sphere
 without first resolving the Kakeya problem.

As might be expected the proof of Theorem 1 uses Kakeya  techniques related to
 Bourgain's paper \cite{B0} and the now classical work of C. Fefferman and Cordoba. 
The  necessary geometric information  while
 not particularly deep is  different from what has been used previously, 
and we prove what we need in section 1 below. In Section 2 we discuss a lemma 
from \cite{M}, in Section 3 we prove our main lemma (Lemma 3.5) and in Section 4 we 
prove Theorem 1. In Section 5 we prove the corollary and make some further 
related remarks. Finally, in an appendix we discuss the related question of 
mixed norm estimates for
the restriction of the X-ray transform to the light rays.
 We prove an optimal local result (except for endpoint questions)
 in three and four dimensions and a partial result in higher dimensions. This is stated below as Theorem A.1.

We will use several  ideas and lemmas from the  previous work on the cone problem, 
e.g.\ from  \cite{B}, \cite{M} and \cite{TV}. Some aspects of the argument and also
 the fact that Theorem 1 should be an accessible result were suggested by the 
author's recent paper \cite{W}.

\medbreak{\it List of notation}.

$Q(N)$: the cube in $\R^d$ centered at the origin with side length $N$.

$|E|$: measure or cardinality of the set $E$ depending on the context.

$\chi_E$: indicator function of $E$.

\section{A property of light rays} 
\advance\eqcount by 1

In this section we fix a suitable large constant $B$ depending on the dimension $d$. 

A {\it light ray} will mean a line in $\R^d$ making a $45$ degree angle 
with the  plane $x_d=0$. We fix two disjoint conical sets $\Gamma_1$ and
 $\Gamma_2$ as described in the introduction and will say that a
light ray is {\it white} (resp.\ black) if its direction belongs to $\Gamma_1$ (resp.\ $\Gamma_2$). 
Thus any white and
 black rays are transverse.  We fix a small positive number  $\e$. 

Let $\d>0$, and let $\W$ and $\B$ be sets respectively of white and black light
 rays with respective cardinalities $m$ and $n$.    For each white  line ${\bf W}$ 
(or black line ${\bf B}$) we  associate to ${\bf W}$ (or ${\bf B}$) the infinite cylinder whose axis is ${\bf W}$ 
(or ${\bf B}$) and whose cross section radius is $\d$. We will denote these tubes 
by $w$ and $b$. For each  tube $w$ (similarly $b$) we define
\begin{equation}\phi_w(x)= \min\left(1,\frac{\d}{\mbox{dist}(x,w)}\right)^M
\label{qq}\end{equation} 
where $M$ is a large constant depending on $\e$. We assume that $\W$  
(similarly $\B$) is {\it $\d$-separated}; by this we  mean the 
following: if $D$ is a disc in projective space with radius $\d$, 
then the tubes $w$ whose axes belong to $D$ 
have bounded overlap, i.e.\ no point belongs to more than $B$ of them. 
We note this implies that the cardinality of lines in
$\W$ which intersect a given compact set is bounded by a (negative)
 power of $\d$.

A {\it $\mu$-fold point} is a point which belongs to at least $\mu$ 
white tubes, and  a {\it smooth $\mu$-fold point} is a point where the quantity 
$$\Phi_{\W}\stackrel{{\rm def}}{=}\sum_{w\in\W}\phi_w$$
is at least equal to $\mu$.

We fix a partition of $Q(1)$ into pairwise disjoint $\d^{\e}$-cubes; in this section
 we  reserve the letter $Q$ for these cubes (except for the standing notation $Q(N)$
 for cubes centered at the origin). In what follows we will be working  with a
 relation $\sim$ between white or black tubes and the cubes $Q$. For any such
 relation we denote
\begin{eqnarray*}
n_{\W}(Q)&=&|\{ w\in\W:w\sim Q\}|,\\
 n_{\B}(Q)&=& |\{ b\in\B:b\sim Q\}|.
\end{eqnarray*} 
 If $x$ is a point or $E$ is a set contained in a cube $Q$ then we will use the
 notation
$$w\sim x \;\; \mbox{ (resp.\ $w\sim E$)}$$
to mean that $w\sim Q$, where $Q$ is the $\d^{\e}$-cube containing $x$ (resp.\ $E$), 
and we define
\begin{eqnarray*}
\tilde{\Phi}_{\W}(x)&=&\sum_{\stackrel{w\in\W}{w\not\sim x}}\phi_w(x),\\
 \tilde{\Phi}_{\B}(x)&=&\sum_{\stackrel{b\in\B}{b\not\sim x}}\phi_b(x).\end{eqnarray*}

We also define (cf. \cite{B0}) a {\it bush}  to be a set of   tubes which are all the 
same color and which all pass through a common point $p$,
and more generally an {\it $\eta$-bush} to be a set of  tubes which are all the same color
 and are all at distance $<\eta$ from a common  point $p$. We call any such point $p$
a {\it base point}
for the bush. 

The purpose of this section is to prove the following lemma.

\nonumproclaim{Lemma 1.1} Assume $\W$ and $\B$ are $\d$\/{\rm -}\/separated{\rm .} Then there is a 
relation $\sim$ between white or black tubes and $\d^{\e}$\/{\rm -}\/cubes $Q$ so that the
following hold{\rm ,} where $C$ depends on $d$ only\/{\rm ;} the implicit constants also depend 
on $\e$\/{\rm :}\/
\medbreak
{\rm 1.}  $\sum_Q n_{\W}(Q)\lesssim m\lo{5}${\rm .}
\smallbreak
{\rm 2.}  $\sum_Q n_{\B}(Q)\lesssim n\lo{5}${\rm .}
\smallbreak\hangindent=38pt\hangafter=1
{\rm 3.} The $\d$-entropy of the set $\{ x\in Q(1):\tilde{\Phi}_{\W}(x)\geq\mu \mbox{ and }\Phi_{\B}(x)\geq\nu\}$ is 
$\lesssim \d^{-C\e}\frac{mn}{\mu^2\nu}${\rm .}
\smallbreak\hangindent=38pt\hangafter=1
{\rm 4.} The $\d$-entropy of the set $\{x\in Q(1):\Phi_{\W}(x)\geq\mu \mbox{ and }\tilde{\Phi}_{\B}(x)\geq\nu\}$ is 
$\lesssim \d^{-C\e}\frac{mn}{\mu\nu^2}${\rm .}

\endproclaim

\demo{{R}emarks} 1. It is easy to see that the $\d$-entropy of the points which
 belong to
$\mu$ white and $\nu$ black tubes can be as large as $\frac{mn}{\mu\nu}$ - just take 
$\W$ and $\B$ to be bushes with a common basepoint and set $\mu=m$, $\nu=n$. Thus property 3 
 gains a factor of $\mu$ over the \lq\lq trivial" bound valid with $\tilde{\Phi}_{\W}$ replaced by $\Phi_{\W}$.
 In the proof of Theorem 1,
this factor will compensate for the factor appearing in  Mockenhaupt's
 estimate for the relevant square function, i.e.\ in Lemma 2.1 below.
 It is also important that the dependence on $\d$ in properties 1 and 2  is only 
logarithmic, or more precisely that it does not involve the specific power $\d^{-\e}$.
On the other hand the distinction between $\mu$-fold points and smooth\break
 $\mu$-fold points
is  purely technical - the functions $\phi_b$ are needed later on in 
order to estimate Schwartz tails.
 \smallbreak
2. It is natural to state Lemma 1.1 in the above manner, since only  
properties 1-4 of the relation $\sim$ will be used in the subsequent sections 
and not its exact definition. However, the relation will be constructed 
in an explicit and  fairly simple  way: roughly, arrange the white or black tubes
 into bushes, and define $w\sim Q$ if 
$w$ belongs to a bush whose basepoint is in $Q$. This procedure together with 
the induction argument in section 4 below
is a variant on the \lq\lq two ends"
argument in \cite{Wr1}, \cite{Wr2}.
\enddemo

Lemma 1.2 below  is  true because $\phi_w$ is essentially a rapidly decreasing sum of
 constants times characteristic functions of dilates of $w$; we leave the details 
to the reader. Lemma 1.3 is a geometrical fact; similar facts are used in 
various places in the literature, e.g.\ in \cite{B} and \cite{TV}.

\nonumproclaim{Lemma 1.2}  If $x\in Q(1)$ is a smooth $\mu$\/{\rm -}\/fold point for the white 
tubes with $\mu\geq\d^B$ then 
 $x$ is a basepoint for 
an $\eta$\/{\rm -}\/bush {\rm (}\/of white tubes\/{\rm )} with cardinality  $\gtrsim
(\log\frac{1}{\d})^{-1}\mu(\frac{\eta}{\d})^{M}$ for some
 $\eta\leq\d^{1-\e}${\rm .} Conversely if $C$ is a large fixed
 constant and $x\in Q(1)$ is a basepoint for an $\eta$\/{\rm -}\/bush 
with cardinality $\geq C\mu(\frac{\eta}{\d})^{M}$  then
$x$ is a smooth $\mu$\/{\rm -}\/fold point{\rm .}
\enddemo

\nonumproclaim{Lemma 1.3} Let $\C\subset\W$ 
be an $\eta$-bush with (say) $\eta\leq\sqrt{\d}$  and  let $p$ be a
 basepoint for $\C${\rm .} Define a set $\Omega$ by deleting from $Q(1)$ the
 double of the $\d^{\e}$\/{\rm -}\/square $Q$ containing $p${\rm .} Let $b$ be any 
black tube{\rm .} Then
\begin{equation}\int_{\Omega}\phi_b\Phi_{\C}\lesssim \d^{-\e(d-2)}\d^d(\frac{\eta}{\d})^{2d-3}.\label{o7}
\end{equation}

\phantom{cliff}

\demo{Proof} First let $b$ and $w$ be a black and a white tube respectively. For any
 $\l\leq 1$ the set 
$$\{x\in Q(1):\phi_b(x)\geq\l\}$$
is contained in a tube with the same axis as $b$ and with width about
 $\d\l^{-\frac{1}{M}}$, and similarly with $w$. Since $w$ and $b$ are
 transverse we have the bound
\begin{equation}|\{x\in Q(1):\min(\phi_b(x),\phi_w(x))\geq\l\}|
\lesssim\left(\d\l^{-\frac{1}{M}}\right)^d\label{qq2} . \end{equation}
 Let $\Delta(b,w)$ be the quantity $\inf_{x\in\Omega}(\mbox{dist}(x, b)
+\mbox{dist}(x, w))+\d$. If $\l$ is large compared with $( {\d}/{\Delta(b,w)})^M$ 
then the set in \bref{qq2} does not intersect $\Omega$. It follows therefore 
that
  \begin{eqnarray}\int_{\Omega}\phi_b\phi_w&\lesssim&
\int_{\Omega}\min(\phi_b,\phi_w)\label{o1}\\ 
&\lesssim&\int_{0}^{({\d}/{\Delta(b,w)})^M}\left(\d\l^{-\frac{1}{M}}\right)^dd\l\nonumber\\&\lesssim&\d^d
\left(\frac{\d}{\Delta(b,w)}\right)^{M-d}.\nonumber
\end{eqnarray}

Now we prove the estimate \bref{o7} when $\eta=\d$. It is clear from \bref{o1}
 that the contribution to the left side from tubes $w\in{\cal C}$ such that $\Delta(b,w)\gtrsim\d^{\e}$ is 
small. On the other hand let $\rho$ be small
 compared with $\d^{\e}$, and consider how many tubes $w\in{\cal C}$ there
 can be with $\Delta(b,w)\leq\rho$. The bush ${\cal C}$
 is clearly contained in a $C\d$-neighborhood of the portion of the light 
cone with origin at $p$ which corresponds to the conical subset $\Gamma_1$. 
If $b$ contains a certain point $y$ which lies within $\rho$ of $\Gamma_1$ and is farther than 
$\d^{\e}$ from $p$,
then by transversality $b$ must intersect  $\Gamma$ at a  point within $C\rho$ of
 $y$. Thus the  number of tubes $w$ with $\Delta(b,w)\leq\rho$ is bounded by the 
$\d$-entropy of the set of lines in $\Gamma$ which intersect a fixed $C\rho$-disc
 lying at distance farther than $\d^{\e}$ from the vertex; equivalently, by the 
$\d$-entropy of a $\d^{-\e}\rho$-disc on $S^{d-2}$, which is 
$({\d^{-\e}\rho}/{\d})^{d-2}$. We conclude
using \bref{o1} that there is a bound
$$\sum\left(\frac{\d^{-\e}\rho}{\d}\right)^{d-2}\d^d\left(\frac{\d}{\rho}\right)^{M-d}$$
with the sum being over dyadic $\rho\geq\d$. Thus we get the bound
$\d^{d-\e(d-2)}$ as claimed.

 We now remove the restriction $\eta=\d$. If $\C$ is an $\eta$-bush then, for
 parameters $\rho$ such that $\rho\geq \eta$ but     small compared with
 $\eta^{\e}$, the maximum number of $\eta$-separated lines in ${\cal C}$ with $\Delta(b,w)\leq\rho$ is bounded by
$({\d^{-\e}\rho}/{\eta})^{d-2}$; for this
 just apply the above argument replacing $\d$ by $\eta$. The space of light\break rays 
is $(2d-3)$-dimensional, so any fixed light ray can be within $\eta$ of at\break most
$(\frac{\eta}{\d})^{2d-3}_{\phantom{|}}$ 
$\d$-separated ones.  It follows that for any 
$\rho \ll \eta^{\e}$ there are\break $\lesssim(\frac{\eta}{\d})^{2d-3}({\d^{-\e}(\rho+\eta)}/{\eta})^{d-2}$  
tubes $w$ with $\Delta(b,w)\leq\rho$. We now 
apply \bref{o1}  as above to bound the left side of \bref{o7} by
$$\sum\left(\frac{\eta}{\d}\right)^{2d-3}\left(\frac{\d^{-\e}(\rho+\eta)}
{\eta}\right)^{d-2}\d^d\left(\frac{\d}{\rho}\right)^{M-d}$$
plus a negligible error, with the sum being over dyadic $\rho\geq\d$.
 Estimate \bref{o7} follows  from this.
\enddemo

The following  lemma is the main step in the argument.
 Essentially, it corresponds to Lemma 1.1 except that here  we ignore the
tails (they will be taken care of in the next lemma) and work with a fixed 
value of $\mu$ (hence the induction argument in the last part of the proof
 of Lemma 1.1 below).

\nonumproclaim{Lemma 1.4} Given a value of $\mu_0$ we can partition $\W$ as
$$\W=\W_g\cup\W_b$$
where 
\begin{itemize}
\ritem{1.} $\W_g$ has no $\mu_0$-fold points in $Q(1)${\rm ,} and

\ritem{2.} $\W_b=\cup_{i=1}^R\C_i$ where each $\C_i$ is a bush with basepoint in $Q(2)$ and\break
$R\lesssim\frac{m}{\mu_0}(\log\frac{1}{\d})^2$. $\phantom{\sum^\int}$
\end{itemize}

\endproclaim

{\it Proof}.
 We fix a large enough constant $C=C_d$ and then another large constant $A$. We will 
use a recursive argument. Accordingly, if $\W^i\subset\W$, then we let $\kappa(\W^i)$ 
be the maximum possible cardinality for a set of $\d$-separated\break $\mu_0$-fold points for
 $\W^i$. We have 
$\kappa(\W)\lesssim \d^{-d}$ since all the tubes in $\W$ are contained in a fixed compact
 set.

Assume now that $\kappa(\W^i)=k$. We will prove: $\W^i=\W^{i+1}\cup\W^i_b$
where $\kappa(\W^{i+1})\leq\frac{k}{2}$, and $\W^i_b$ is the union of $\lesssim A\frac{m}{\mu_0}
\log \frac{1}{\d}$ $\d$-bushes.

Namely, let ${\cal R}_i$ be a set of $\d$-separated $\mu_0$-fold points for $\W_i$ with 
maximum possible cardinality $k$. There are two cases.
\begin{itemize}
\item[(i)] If $k\leq A\frac{m}{\mu_0}\log \frac{1}{\d}$ then we let $\W^i_b$ be all tubes 
 $w\in\W^i$ such that\break dist$(x,w)<\d$ for some $x\in {\cal R}_{i}$ and $\W^{i+1}=\W^i\backslash\W^i_b$. 
Evidently $\W^i_b$ is the union of $\lesssim A\frac{m}{\mu_0}\log \frac{1}{\d}$  $\d$-bushes; and 
$\kappa(\W^{i+1})=0$ since 
any $\mu$-fold point for $\W^i$ must lie within $\d$ of some 
point of  ${\cal R}_i$.

\item[(ii)] If $k> A\frac{m}{\mu_0}\log \frac{1}{\d}$ we choose $A\frac{m}{\mu_0}\log 
\frac{1}{\d}$ points  from ${\cal R}_i$ at random. We let $\W_b^i$ be the tubes 
 $w\in\W^i$ such that dist$(x,w)<\d$ for some $x$  in the random sample, and $\W^{i+1}=\W^i\backslash\W^i_b$.
 Evidently $\W^i_b$ is the union\break of $\lesssim A\frac{m}{\mu_0}\log \frac{1}{\d}$ $\d$-bushes. We will show that 
with high 
probability\break $\kappa(\W^{i+1})\leq \frac{k}{2}$.
\end{itemize}

For this, define for each $w\in\W^i$
$$P(w)=k^{-1}|\{x\in{\cal R}_i: \mbox{dist}(w, x)<\d\}|.$$
 Thus the probability that $w$ is in $\W^{i+1}$ is at most
$$(1-P(w))^{A\frac{m}{\mu_0}\log \frac{1}{\d}}.$$ 
If $P(w)\geq C^{-1}\frac{\mu_0}{m}$ it follows that the probability that $w$ is
 in $\W^{i+1}$ is at most
 $\d^{\frac{A}{C}}$. If $A$ is large enough then since the cardinality of the set 
of lines in $\W$ which intersect $Q(2)$ is bounded by $\d^{-B}$ it follows that 
with high probability no tubes with $P(w)\geq C^{-1}\frac{\mu_0}{m}$ belong to
$\W^{i+1}$. 

Now let ${\cal R}_{i+1}$ be a maximal set of $\d$-separated $\mu_0$-fold points 
for $\W^{i+1}$, and let ${\cal R}$ be a maximal $2\d$-separated subset of ${\cal R}_{i+1}$. Consider
 the  quantity
\begin{equation}\sum_{w\in\W^{i+1}}kP(w)\label{cl1}.
\end{equation}
We have seen that  with high probability \bref{cl1} is less than $\frac{k}{C}\frac{\mu_0}{m}|\W^{i+1}|
\leq\frac{k}{C}\mu_0$. On the other hand,
 we have
\begin{eqnarray*}
\noalign{\vskip9pt} \sum_{w\in\W^{i+1}}kP(w) &=&\sum_{x\in{\cal R}_i}|\{w\in \W^{i+1}:\mbox{dist}
(w, x)<\d\}|\\ \noalign{\vskip9pt}
&\geq&\mu_0|{\cal R}|.\\ \noalign{\vskip-6pt}
\end{eqnarray*} 
The first line follows from the 
definition   by reversing the order of summation, and the second line 
then follows because every point in ${\cal R}$ is within $\d$ of a point of
 ${\cal R}_i$ and no two points of ${\cal R}$ can be within $\d$ of the same point 
of ${\cal R}_i$. We conclude that with high probability $|{\cal R}|\leq C^{-1}k$. Since 
$|{\cal R}|$ and $|{\cal R}_i|$ are comparable it then follows that $|{\cal R}_{i+1}|\leq\frac{k}{2}$, 
as was to be shown.

We now proceed recursively. Let $\W^0=\W$ and apply the preceding to express $\W^0=\W^0_b\cup\W^1$.
 Then apply the 
preceding
to express $\W^1=\W^1_b\cup\W^2$ and continue in this manner, stopping when we reach a
 situation where we are in 
case (i) above. Suppose we stop after $T$ stages. Since $\kappa(\W^i)$ is initially 
$\lesssim \d^{-d}$ and decreases each time 
at least by a factor of $2$, we then have $T\lesssim\log \frac{1}{\d}$. We now define
 $\W_g$ to be the set $\W^{i+1}$ defined 
at the last iteration. It satisfies $\kappa(\W_g)=0$ as required. On the other hand we
 define $\W_b=\cup_{i}\W^i_b$.
 This set is the union of the $\lesssim\log \frac{1}{\d}$ sets $\W^i_b$, each of which  is 
the union of 
$\lesssim A\frac{m}{\mu_0}\log \frac{1}{\d}$ bushes. The lemma follows.\hfill\qed
 
 \phantom{wet}


The next lemma is a version of the preceding one incorporating  Schwartz tails. 

 \nonumproclaim{Lemma 1.5} Fix $\mu_0\geq\d^{B}$. Then $\W=\W_g\cup\W_b$ where
\begin{itemize}
\ritem{1.} $\Phi_{\W_g}\leq\mu_0$ everywhere{\rm ,}

\ritem{2.} $\W_b=\cup_{k:2^k<\d^{-\e}}\W_b^k${\rm ,} and for each $k${\rm ,} $\W_b^k=\cup_{i=1}^{R_k}\C_i${\rm ,}
 where $\C_i$ is a $2^k\d$\/{\rm -}\/bush with basepoint in $Q(2)$ and $R_k\lesssim\frac{m}{2^{Mk}\mu_0}\lo{3}${\rm .}
\end{itemize}

\endproclaim

\demo{Proof}  Let $w^{\eta}$ be the $\eta$-tube with the same axis as $w$. Notice
 that Lemma 1.3 is applicable also to
the $w^{\eta}$'s (provided $\log\eta$ is comparable to $\log\d$, which will be the case 
below), since we used $\d$-separation in the proof only to conclude that the cardinality
 of the white lines which intersect $Q(2)$ was bounded by a negative power of $\d$.

We now define recursively a family of subsets $\W_b^j$. Let $\W=\W_g^1\cup\W_b^1$ be
 the decomposition from Lemma 1.3 for the given $\mu_0$. If $k\geq 2$ and if $\W_b^j$
 have been defined for $j<k$ then we let $\W_g^{k-1}=\W\backslash(\cup_{j=1}^{k-1}\W_b^j)$.
 The following inductive hypothesis will hold:
\bigbreak
$(\ast)$  If $j\leq k-1$, then the family of tubes $\{w^{2^j\d}: w\in \W_g^{k-1}\}$ has 
no $2^{Mj} \mu_0$-fold points.
\bigbreak

 Let $\eta=2^k\d$ and apply Lemma 1.3 to
the tubes $\{w^{\eta}:w\in\W_g^{k-1}\}$ replacing $\mu_0$ by $2^{Mk}\mu_0$.  This 
decomposes $\W_g^{k-1}=\W_b^k\cup\W_g^k$ where the tubes\break $\{w^{2^k\d}:w\in\W_g^k\}$ 
have no
$2^{Mk}\mu_0$-fold points and $\W_b^k$ is the union of at most $\frac{m}{2^{Mk}\mu_0}
\lo{2}$ $2^{k}\d$-bushes. The inductive hypothesis $(\ast)$ is then satisfied for 
$j\leq k$. We continue in this manner, stopping when $2^{Mk}\mu_0$ becomes greater 
than $m$. This will occur at a stage $k$ with $2^k<\d^{-\e}$, since we have assumed $\mu_0\geq\d^B$. We
 define $W_g$ to be the last $\W_g^k$.

If $\Phi_{\W_g}(x)\geq (C\log\frac{1}{\d})\mu_0$ with $C$ a large fixed constant then by Lemma~1.2 $x$ must be a
$2^{Mk}\mu_0$-fold point for the tubes $\{w^{2^k\d}:w\in\W_g\}$ for some $k$, hence also a $2^{Mk}\mu_0$-fold point for
the larger family $\{w^{2^k\d}:w\in\W_g^k\}$, which is impossible by construction.  The lemma now follows by replacing
$\mu_0$ with $(C\log\frac{1}{\d})^{-1}\mu_0$. \enddemo

 \vglue9pt

To prove Lemma 1.1 it suffices by symmetry to construct a relation between white tubes and $\d^{\e}$-squares 
so that  properties 1  and 3  hold. This will again be done recursively. A remark on terminology: in this argument,\break
when we  say  that \lq\lq$\C$ is a $2^k\d$-bush"
we mean that $\C$ is a $2^k\d$-bush but not a $2^{k-1}\d$-bush.

We apply Lemma 1.5 to $\W$ with $\mu_0=\frac{m}{2}$, obtaining a set $\W_g^1$
with\break $\Phi_{\W_g^1}\leq\frac{m}{2}$  and 
a collection of stage 1 $\eta$-bushes $\C_i^1$ (thus each $\C_i^1$ is a\break $2^k\d$-bush for some $k$ with $2^k\d\leq
\d^{1-\e}$). Then we apply Lemma 1.5 to $\W_g^1$ with $\mu_0=\frac{m}{4}$ obtaining $\W_g^2$ with $\Phi_{\W_g^2}\leq\frac{m}{4}$ and stage two $\eta$-bushes $\C_i^2$ and continue in this manner, taking 
$\mu_0=\frac{m}{2^j}$ at the $j$th stage. We stop the induction at stage $R$, where $R$ by definition is the 
smallest integer such that $\frac{m}{2^R}<\d^B$. Clearly $R\lesssim\log\frac{1}{\d}$.  For each $j_0\leq R$ we 
now have a decomposition \vglue-9pt
\begin{equation}\W=\W_g^{j_0}\cup(\cup_{j\leq j_0}\cup_i\C_i^j)\label{18.1}\end{equation}
\vglue6pt\noindent
where $\Phi_{\W_g^{j_0}}\leq\frac{m}{2^{j_0}}$, and (by part 2  of Lemma 1.5) we have the following:
\begin{itemize}
 
\item[] For each  $j$ and $k$ there are $\lesssim 2^j2^{-Mk}\lo{3}$ values of $i$ such that $\C_i^j$ is a $2^k\d$-bush.
\item[]
 For each $\C_i^j$ we fix a basepoint $p_i^j$. We now define the relation $\sim$:
\end{itemize}

\demo{Definition} A tube $w$ and $\d^{\e}$-square $Q$ are related, $w\sim Q$,  if $w$ belongs to an 
$\eta$-bush $\C_i^j$ such that $p_i^j$ is in $Q$ or one of its neighbors. 
\enddemo

We show first that  property 1  holds. Suppose that $\C_i^j\subset\W_g^{j-1}$ is a\break $2^k\d$-bush. Then, using Lemma
1.2  and the fact that $\Phi_{\W_g^{j-1}}\leq\frac{m}{2^{j-1}}$, we get the following bound for the cardinality of 
$\C_i^j$: 
\vglue-5pt
$$|\C_i^j|\lesssim 2^{Mk}\Phi_{\C_i^j}(p_i^j)\lesssim 2^{Mk}\frac{m}{2^j}.$$
\vglue6pt\noindent
By the preceding bound for the number of $2^k\d$-bushes, we then have $$\sum_i|\C_i^j|\lesssim \sum_k2^j2^{-Mk}
\lo{3}\cdot  2^{Mk}\frac{m}{2^j}\lesssim m\lo{4}.$$
 Summing over $j$ we get $\sum_{i,j}|\C_i^j|\lesssim m\lo{5}$. Thus, there are at most $m\lo{5}$ pairs $(w,\C)$ where
 $w$ is a white tube and $\C=\C_i^j$ is an $\eta$-bush containing $w$. This obviously implies property 1. It remains to 
prove property~3. 

Fix $\mu$. If $\mu\lesssim \d^B$ (and if $B=B_d$ was chosen large enough) then property~3  will clearly hold, since the
right
 hand side will be greater than $\d^{-d}$. On the other hand, if $\mu$ is large compared with $\d^B$ then we can
 choose $j_0$ so that $\frac{m}{2^{j_0}}$ is less than $\frac{\mu}{2}$ but greater than $\frac{\mu}{8}$. We consider
 the decomposition  \bref{18.1} with this value of $j_0$. Thus
 $\Phi_{\W_g^{j_0}}\leq\frac{\mu}{2}$ and  for each $k$ we have
\begin{equation}|\{(i,j):j\leq j_0 \mbox{ and }\C_i^j \mbox{ is a } 2^k\d\mbox{-bush}\}|\lesssim
2^{-Mk}\lo{3}\frac{m}{\mu}.
\label{18.2}\end{equation} 
Fix a black 
tube $b$, and  fix also a choice of  $\C_i^j$ with $j\leq j_0$. Define $\Omega_{ij}$ by deleting from $Q(1)$ the
 $\d^{\e}$-square 
containing $p_i^j$ and its neighbors. Lemma 
1.3 implies that if $\C_i^j$ is a $2^k\d$-bush then
$$\int_{\Omega_{ij}}\phi_b\Phi_{\C_i^j}\lesssim \d^{-(d-2)\e}2^{Ck}\d^d$$
where $C$ depends on $d$.

Now sum over $b$, $i$ and $j\leq j_0$ obtaining (provided $M$ has been chosen large enough)
\begin{eqnarray*}
\sum_{ij}\int_{\Omega_{ij}}\Phi_{\B}\Phi_{\C_i^j}&\lesssim& n\sum_k2^{-Mk}\lo{3}\frac{m}{\mu}\cdot
\d^{-(d-2)\e}2^{Ck}\d^d\\
&\lesssim &\d^{-(d-2)\e}\frac{nm}{\mu}\d^d\lo{3}
\end{eqnarray*}
where the first inequality follows from
\bref{18.2}.

Suppose now that $ x$ is a point such that $\tilde{\Phi}_{\W}(x)\geq\mu$. By the definition of the relation $\sim$
we have 
$$\tilde{\Phi}_{\W}(x)\leq\Phi_{\W_g^{j_0}}(x)+\sum_{\stackrel{j\leq j_0}{x\in\Omega_{ij}}}\Phi_{\C^j_i}(x).$$
The first term on the right side is $\leq\frac{\mu}{2}$, so 
$$\tilde{\Phi}_{\W}(x)\leq 2\sum_{\stackrel{j\leq j_0}{x\in\Omega_{ij}}}
\Phi_{\C^j_i}(x)$$
whence
$$\int\Phi_{\B}\tilde{\Phi}\leq
2\sum_{ij}\int_{\Omega_{ij}}\Phi_{\B}\Phi_{\C_i^j}\lesssim\d^{-(d-2)\e}\frac{nm}{\mu}\d^d\lo{3}.
$$
 It follows that the measure of  the set where $\Phi_{\B}\geq\nu$ and 
$\tilde{\Phi}\geq\mu$ is $$\lesssim\d^{-(d-2)\e}\frac{nm}{\nu\mu^2}\d^d\lo{3}.$$

Using that the functions $\phi_w$ are roughly constant on $\d$-discs it 
then follows that the $\d$-entropy is $$\lesssim\d^{-(d-2)\e}\frac{nm}{\nu\mu^2}\lo{3}$$
 as claimed. 
 \enddemo

What we actually use below is a slight variant on Lemma 1.1 where the infinite cylinders 
are replaced by finite ones.
We introduce the following notation which will also be used in Section 3.

\demo{Definition} 1. Suppose that $g$ is a radial function in $\R^d$ and $R$ is 
a centered compact  convex set. Then we use the notation $g_R$ to mean $g\circ A$, 
where $A$ is an affine function mapping (the John ellipsoid for) $R$ onto the unit ball.

2. $\phi$ will denote the function $\phi(x)=\min(1, |x|^{-M})$, where $M$ is a 
sufficiently large constant.
\enddemo

Suppose now that we have   collections $\B$ and $\W$ of cylinders of length $1$ 
and cross section radius $\d$, which are $\d$-separated in the same sense as before;
 i.e.\ the ones whose direction belongs to a given $\d$-disc in projective space have
 bounded overlap, and furthermore the axis directions belong to $\Gamma_1$ and 
$\Gamma_2$ respectively. Let $m=|\W|$, $n=|\B|$. Fix 
(in addition to $\e$) another small positive $\eta$; the choice of $M$ and the 
implicit constants below may now also depend on $\eta$. The quantities $\Phi_{\W}$ and $\tilde{\Phi}_{\W}$
are  defined in the same way as before, except of course that we use the modified definition
of $\phi_w$ via the definition above. 
 
\nonumproclaim{Lemma 1.1$'$}
With the above assumptions there is a relation $\sim$ between white or black tubes 
$w\in\W$ or $b\in\B$ and $\d^{\e}$\/{\rm -}\/cubes $Q\subset Q(1)$ so that the
following hold{\rm ,} where $n_{\W}(Q)=|\{w:w\sim Q\}|$\/{\rm :}
\begin{itemize}
\ritem{1.}  $\sum_Q n_{\W}(Q)\lesssim m\d^{-\eta}${\rm .}

\ritem{2.}  $\sum_Q n_{\B}(Q)\lesssim n\d^{-\eta}${\rm .}

\ritem{3.} The $\d$-entropy of the  set $\{ x\in Q(1):\tilde{\Phi}_{\W}(x)\geq\mu \mbox{ and }\Phi_{\B}(x)\geq\nu\}$ 
is 
$\lesssim \d^{-C\e}\frac{mn}{\mu^2\nu}${\rm .}

\ritem{4.} The $\d$-entropy of the set $\{x\in Q(1):\Phi_{\W}(x)\geq\mu \mbox{ and }\tilde{\Phi}_{\B}(x)\geq\nu\}$ is 
$\lesssim \d^{-C\e}\frac{mn}{\mu\nu^2}${\rm .}
\end{itemize}

\endproclaim

To prove this we define $w\sim Q$ if the infinite cylinder\footnote{We allow the possibility that an infinite cylinder
 may contain several $w$'s. It is therefore easy to reduce to the case where the infinite cylinders are 
$\d$-separated.} with the same axis as $w$ is related to
 $Q$ in the sense of Lemma 1.1 and if in addition the distance from $w$ to the origin
 is less than $\d^{-\frac{\eta}{2}}$. Then properties 1  and 2   follow immediately from properties 1  and 2  of Lemma 1.1, 
 and properties 
3  and 4  follow from properties 3  and 4  of Lemma 1.1 using that the contribution to $\Phi$ from tubes 
further than $\d^{-\frac{\eta}{2}}$ from the origin is negligibly small if $M$ 
is large. \phantom{howtodo} \hfill\qed

\section{A lemma of Mockenhaupt}
\advance\eqcount by 8

 We cover the unit sphere $S^{d-2}$ with a family of spherical caps $c$ of radius 
$N^{-\frac{1}{2}}$ with bounded overlap; this gives also a covering  of $\Gamma$ by 
a family of \lq\lq sectors" $\rho=\rho_c$, where $\rho_c=\{ x\in \Gamma:\frac{\ov{x}}
{x_d}\in c\}$.

We will be using a variant on the square function estimate in \cite{M}. To state it, 
let $\{\rho_j\}$ be the sectors $\rho$ which intersect
$\Gamma_1$ and let $\{\tilde{\rho}_k\}$ be the sectors which intersect $\Gamma_2$. Let 
$f$ and $g$ be two functions on $\R^d$ and assume that $f=\sum_{j=1}^{\mu}f_j$
 and $g=\sum_{k=1}^{\nu} g_k$, where supp$f_j$ is contained in the 
$N^{-{1}/{2}}$-neighborhood of the sector  $\rho=\rho_j$, and
 likewise supp$g_k$ is contained in the $N^{-{1}/{2}}$-neighborhood of
   $\tilde{\rho}_k$. Let $F=\hat{f}$, $G=\hat{g}$, and  
$SF=(\sum_j|\hat{f_j}|^2)^{\frac{1}{2}}$ $SG=(\sum_k|\hat{g_k}|^2)^{\frac{1}{2}}$.

\nonumproclaim{Lemma 2.1} $\|FG\|_2^2\lesssim\min(\mu,\nu)\|(SF)(SG)\|_2^2${\rm .}
\endproclaim 

\demo{Proof} \cite{M} We claim that for a given point $z\in\R^d$ there
 are 
$\lesssim\min(\mu,\nu)$ pairs $(j,k)$ such that $z\in\su f_j+\su g_k$.

We will use the following geometrically obvious fact (a consequence of the 
strict convexity of the sphere): let $\e_0$ be a fixed positive constant and
 let $\zeta, \omega_1,\omega_2$ be points of $S^{d-2}$ with $|\omega_i-\zeta|
\geq\e_0$ for $i=1,2$.
Let $\ell$ be a line  in $\R^{d-1}$ which passes through the point $\zeta$ and 
assume that both $\omega_1$ and $\omega_2$ are at distance at most $\d$ from 
$\ell$. Then $|\omega_1-\omega_2|\leq C\d$, where $C$ depends on $\e_0$.

In order to prove the claim it suffices to show that for fixed $j$ the set of $k$ 
such that $z\in\su f_j+\su g_k$ has bounded cardinality. To this end we fix $\zeta$ 
with $(\zeta,1)\in\rho_j$, and $\omega_1$ and $\omega_2$ such that
$(\omega_i,1)\in \tilde{\rho}_k$ and $z\in\su f_j+\su g_{k_i}$ for $i=1,2$. If we
 let $z=(w, t)$ then for suitable $a,b\in[\frac{1}{2},2]$ we have
$$ a+b=t+{\cal O}(N^{-\frac{1}{2}})$$
$$a\omega_1+b\zeta=w+{\cal O}(N^{-\frac{1}{2}})$$
and therefore
$$a\omega_1+(t-a)\zeta=w+{\cal O}(N^{-\frac{1}{2}})$$
so that
\begin{equation}\omega_1-\zeta=a^{-1}(w-t\zeta)+{\cal O}(N^{-\frac{1}{2}}) . \label{Muck}\end{equation}
Estimate \bref{Muck} says that the distance from $\omega_1$ to the line through
 $\zeta$ spanned by $w-t\zeta$ is $\lesssim N^{-\frac{1}{2}}$. Likewise the 
distance from $\omega_2$ to this line is $\lesssim N^{-\frac{1}{2}}$. The 
disjoint conical support assumption  implies that $|\omega_i-\zeta|$ is bounded
 below for each $i$ so we conclude that $|\omega_1-\omega_2|\leq CN^{-\frac{1}{2}}$.
 This means that there are at most a bounded number of possible values for $k$, 
proving the claim.

The claim implies the lemma by a well-known calculation with the Plancherel theorem,
which we omit.\enddemo

\section{Main lemma} 
\advance\eqcount by 9

 It will  be convenient to change the setup described in the introduction slightly
 in this section. We fix a scale $N$,
 let  $Q(N)$ be the square centered at the origin  with side $N$, and let $\Gamma^{(N)}$
 be the $\frac{1}{N}$-neighborhood
 of $\Gamma_1$; similarly $\Gamma_1^{(N)}$ is the $\frac{1}{N}$-neighborhood of 
$\Gamma_1$, etc. Corresponding to
 the covering of $\Gamma$ by sectors described in Section 2 is a covering  of  
$\Gamma^{(N)}$ by  $\frac{1}{N}$-neighborhoods of sectors, and in this section 
we use $\rho$ to denote one of the latter.
 Thus $\rho$ is essentially a $1 \times
\overbrace{N^{-\frac{1}{2}}\times  \ldots\times N^{-\frac{1}{2}}}^{d-2\, times}
\times N^{-1}$-rectangle. We fix 
disjoint sets $E_{\rho}\subset\rho$ with $\cup_{\rho}E_{\rho}=\Gamma^{(N)}$ 
and let $\zeta_{\rho}=\chi_{E_{\rho}}$.

 Let $f$
 be a function supported on  $\Gamma_1^{(N)}$ with $L^2$ norm $1$, $F=\hat{f}$, $F_{\rho}=\widehat{\zeta_{\rho}f}$, 
and
$$SF(x)=\left(\sum_{\rho}|F_{\rho}(x)|^2\right)^{\frac{1}{2}}.$$

Further let $b$ be a fixed radial Schwartz function nonzero on $Q(1)$ whose Fourier transform has compact
 support and whose $\Z^d$ translations form a partition of unity. For each $\rho$
 we fix a tiling ${\cal F}^{\rho}$ of $\R^d$ by
rectangles  $\s$ with dimensions $N\times\overbrace{N^{\frac{1}{2}} \times 
\ldots\times N^{\frac{1}{2}}}^{d-1\, times}$,
the long direction being orthogonal to the light cone $\Gamma$ at points of 
(the center line of) $\rho$, and we let  $\F=\cup_{\rho}\F^{\rho}$. We also 
let ${\cal P}^{\rho}$ be a tiling by  $N\times \overbrace{N^{\frac{1}{2}} 
\times \ldots\times N^{\frac{1}{2}}}^{d-2\, times}\times 1$ rectangles dual to the sector
 $\rho$. For each $\rho$ and each $\s\in\F^{\rho}$ we define $F_{\rho}^{\s}
=b_{\s}F_{\rho}$, where $b_{\s}$ (and also $b_{\pi}$, $\phi_{\s}$, etc. in 
the subsequent argument) are as in the definition at the end of Section 1; 
thus $\sum_{\s}F_{\rho}^{\s}=F_{\rho}$. For each $(\rho,\s)$ we also further
 decompose $F_{\rho}^{\s}$ as $\sum_{\pi\in{\cal P}^{\rho}}F_{\rho}^{\s,\pi}$,
 where $F_{\rho}^{\s,\pi}=b_{\pi}F_{\rho}^{\s}$. The following fact (trivial to
 prove, since $b$ has compact support) will be very important below:

\nonumproclaim{Lemma 3.1} The inverse Fourier transforms of the functions $F_{\rho}^{\s}$ 
and $F_{\rho}^{\s,\pi}$ are supported in a fixed dilate $\ov{\rho}$ of $\rho${\rm ,} and
 in particular are supported in the $\frac{C}{N}$\/{\rm -}\/neighborhood of $\Gamma${\rm .}
\endproclaim

The following fact is also clear from the Schwartz inequality since $\sum_{\s\in{\cal F}^{\rho}}\phi_{\s}^2$
  and $\sum_{\pi\in{\cal P}^{\rho}}\phi_{\pi}^2$ are bounded for fixed $\rho$. Suppose that for each $\rho$ 
a subset ${\cal A}^{\rho}\subset{\cal F}^{\rho}$ is given. Then
\begin{equation}\sum_{\rho} \bigg|\sum_{\s\in{\cal A}^{\rho}}F_{\rho}^{\s} \bigg|^2
\lesssim \sum_{\rho}\sum_{\s\in{\cal A}^{\rho}}|F_{\rho}^{\s}|^2\phi_{\s}^{-2}
\lesssim\sum_{\stackrel{\rho,\s,\pi}{\stackrel{\s\in{\cal A}^{\rho}}{\pi\in{\cal
P}^{\rho}}}}|F_{\rho}^{\s,\pi}|^2\phi_{\pi}^{-2}\phi_{\s}^{-2}.\label{u1}\end{equation}
\vglue-2pt

The next two lemmas keep track of some relationships among the various 
decompositions of $F$ which follow from orthogonality considerations and the
 uncertainty principle. We note the following:  let $\pi_0$ be a rectangle
 containing the origin, and let $\pi$ be a translate of $\pi_0$. Then,  the
 operator with kernel
\smallbreak
\centerline{${\displaystyle K(x,y)=\phi_{\pi}(x)^{-2}\phi_{\pi_0}(x-y)^{100}\phi_{\pi}(y)^4}$}
\smallbreak\noindent
maps $L^2$ to $L^{\infty}$ with norm $\lesssim |\pi|^{\frac{1}{2}}$, since
 one can easily show that $\int|K(x,y)|^2dy\break\lesssim |\pi|$ for fixed $x$.

\nonumproclaim{Lemma 3.2}  For fixed $\rho$ and $\s\in\F^{\rho}$ we have
 $$\sum_{\pi}\|\phi_{\pi}^{-2}\phi_{\s}^{-3}F_{\rho}^{\s,\pi}\|^2_{\infty}
\lesssim N^{-\frac{d}{2}}\|\phi_{\s}^{-4}F_{\rho}^{\s}\|_2^2.$$
\endproclaim
\vglue-8pt
{\it Proof}. Fix a Schwartz
 function $\kappa$ whose Fourier transform is $1$ on the unit ball and let
 $\ov{\kappa}$ be the corresponding 
function whose Fourier transform is $1$ on  the set $\ov{\rho}$ in Lemma 
3.1, obtained from $\kappa$ by composition with a linear map followed by
 multiplication by a character and by a scalar with magnitude about $|\rho|$.
 Then 
$F_{\rho}^{\s,\pi}=\ov{\kappa}\ast F_{\rho}^{\s,\pi}$. Let $\s_0$  and 
 $\pi_0$  be the  rectangles in the tilings ${\cal F}^{\rho}$ and ${\cal P}^{\rho}$ 
which contain the origin. Then $|\ov{\kappa}(z)|\lesssim
 |\rho|\phi_{\pi_0}(z)^{200}\lesssim|\rho|\phi_{\pi_0}(z)^{100}\phi_{\s_0}(z)^{100} $.
 We conclude that
\smallbreak
\centerline{${\displaystyle |\phi_{\pi}^{-2}(x)\phi_{\s}^{-3}(x)F_{\rho}^{\s,\pi}(x)|\lesssim \int K(x,y)|
\phi_{\pi}^{-4}(y)\phi_{\s}^{-4}(y)F_{\rho}^{\s,\pi}(y)|dy}$}
\smallbreak\noindent
where
\begin{eqnarray*}K(x,y)&=&|\rho|\phi_{\pi}(x)^{-2}\phi_{\s}(x)^{-3}
\phi_{\pi_0}(x-y)^{100}\phi_{\s_0}(x-y)^{100}
\phi_{\pi}(y)^4\phi_{\s}(y)^4\\ \noalign{\vskip6pt}
&\lesssim&|\rho|\phi_{\pi}(x)^{-2}\phi_{\pi_0}(x-y)^{100}\phi_{\pi}(y)^{4}.
\end{eqnarray*}
We have seen that the norm of this kernel from $L^2$ to $L^{\infty}$ is 
$\lesssim |\pi|^{\frac{1}{2}}|\rho|\approx|\rho|^{\frac{1}{2}}\approx 
N^{-\frac{d}{4}}$.
Accordingly
 \begin{eqnarray*}\sum_{\pi}\|\phi_{\pi}^{-2}\phi_{\s}^{-3}
F_{\rho}^{\s,\pi}\|^2_{\infty}&\lesssim &N^{-\frac{d}{2}}
\sum_{\pi}\|\phi_{\pi}^{-4}\phi_{\s}^{-4}F_{\rho}^{\s,\pi}\|_2^2\\
&=&N^{-\frac{d}{2}}
\sum_{\pi}\|(\phi_{\pi}^{-4}b_{\pi})\phi_{\s}^{-4}F_{\rho}^{\s}\|_2^2\\
\noalign{\noindent and now we use that $\sum_{\pi}|\phi_{\pi}^{-4}b_{\pi}|^2\lesssim 1$ 
pointwise, 
obtaining the lemma.\hfill\qed}
\end{eqnarray*}
  
\phantom{strange}
\vglue-36pt
For each $\rho$ and each $\s\in\F^{\rho}$  we define a parameter
$$ h(\s)=(N^{-\frac{d+1}{2}}\|\phi_{\s}^{-4}F_{\rho}^{\s}\|_2^2)^{\frac{1}{2}}.$$
We think of $h(\s)$ as being essentially the $L^2$ average of $F_{\rho}$ on 
$\s$.
We group the $\s$'s into families corresponding to the different possible dyadic 
values for $h(\s)$; thus
$$\F(h)=\{\s\in \F:h(\s)\in [\frac{h}{2}, h]\}$$
and we define
$$F_h=\sum_{\rho}\sum_{\s\in\F(h)\cap\F^{\rho}}F_{\rho}^{\s}.$$
 \nonumproclaim{Lemma 3.3}  $h^2|\F(h)|\lesssim N^{-\frac{d+1}{2}}${\rm .}
\endproclaim
 
\demo{Proof} Clearly
$$h^2|\F(h)|\lesssim N^{-\frac{d+1}{2}}\sum_{\rho}\sum_{\s\in\F^{\rho}}\|\phi_{\s}^{-4}F_{\rho}^{\s}\|_2^2.$$
For fixed $\rho$ we have $\sum_{\s}|b_{\s}\phi_{\s}^{-4}|^2\lesssim 1$ pointwise.
 So for fixed $\rho$ we have $\sum_{\s\in\F^{\rho}}\|\phi_{\s}^{-4}
F_{\rho}^{\s}\|_2^2\lesssim \|F_{\rho}\|_2^2$. If we sum over $\rho$ and
 use orthogonality of the $F_{\rho}$'s
 the lemma follows. \enddemo

 If  $g$ is a function supported on $\Gamma_2^{(N)}$ with $L^2$ norm $1$ we will 
likewise denote $\widehat{gd\s}$ by $G$, etc. Thus we obtain also functions $G_{\rho}$, $G_{\rho}^{\s}$,
 $G_{\rho}^{\s,\pi}, G_h$, and families of tubes $\G$, $\G^{\rho}$, $\G(h)$.
 The next lemma is a  \lq\lq local" estimate; it will then be combined with Lemma 1.1 
to give the following Lemma 3.5 which is the main result of this section. 
 
\nonumproclaim{Lemma 3.4} Fix  a square $Q$ with side $\sqrt{N}${\rm .} Let $\tilde{\F}$ and 
$\tilde{\G}$ be subsets of $\F(h_1)$ and $\G(h_2)$ respectively and let $\mu$ and 
$\nu$ be the maximum values on the square $Q$ of the functions $\Phi_{\tilde{\F}}$
 and $\Phi_{\tilde{\G}}${\rm .} Then
\begin{equation}\int_Q|(\sum_{\rho}\sum_{\s\in\tilde{\F}\cap\F^{\rho}} F_{\rho}^{\s})(\sum_{\rho_2}\sum_{\s_2\in\tilde{\G}\cap\G^{\rho_2}}G_{\rho_2}^{\s_2})|^2
\lesssim h_1^2h_2^2\mu\nu\min(\mu,\nu)
N^{\frac{d}{2}}.\qquad\label{2.1.1}\end{equation}
\endproclaim

\demo{Proof} 
We subdivide $\tilde{\F}$ and $\tilde{\G}$ according to the possible dyadic values 
for $\phi_{\s}$ on 
$Q$. Thus we define 
$$\tilde{\F}(k)=\{\s\in\tilde{\F}:\min_Q\phi_{\s}\in[2^{-(k+1)},2^{-k}]\}$$
$$\tilde{\G}(\ell)=\{\s\in\tilde{\G}:\min_Q\phi_{\s}\in[2^{-(\ell+1)},2^{-\ell}]\}.$$
 We note that if
$\s\in\tilde{\F}(k)$ then \begin{equation}\|\phi_{\s}\phi_Q\|_{\infty}\lesssim 2^{-k}.\label{bref}\end{equation} 
This follows from the rapid decay of $\phi$ and the fact that $\s$ contains a translate 
of $Q$. Hence also $\|\phi_{\s}b_Q\|_{\infty}\lesssim 2^{-k}$. Furthermore, from the 
definition of $\mu$ and $\nu$, we have 
\begin{equation} |\tilde{\F}(k)|\lesssim 2^k\mu \mbox{ and }|\tilde{\G}(\ell)|\lesssim 2^{\ell}\nu.\label{19.1}
\end{equation}

The left side of \bref{2.1.1} is $\lesssim \sum_{k=0}^{\infty}\sum_{\ell=0}^{\infty} 
2^{k+\ell}A(k,\ell)$, where
\begin{equation}A(k,\ell)=\int\left|\left(b_Q^3\sum_{\rho}\sum_{\s\in\tilde{\F}(k)\cap 
\F^{\rho}} F_{\rho}^{\s}\right)\left(b_Q^3\sum_{\rho_2}\sum_{\s_2\in\tilde{\G}(\ell)\cap\G^{\rho_2}}
G_{\rho_2}^{\s_2}\right)\right|^2. \qquad\label{pro5}\end{equation}
 Using Lemma 3.1 and that $\hat{b}$ has compact support, one sees that the Fourier 
transform of the function $b_Q^3\sum_{\s\in\tilde{\F}(k)\cap \F^{\rho}} F_{\rho}^{\s}$ 
is supported in the $CN^{-\frac{1}{2}}$-neigh\-borhood of the sector $\rho$; and similarly
 with the second factor in \bref{pro5}. Lemma~2.1 is therefore applicable and implies
 that
$$A(k,\ell)\lesssim
\min(|\tilde{\F}(k)|,|\tilde{\G}(\ell)|)\int\sum_{\rho}\left|b_Q^3\sum_{\s\in\tilde{\F}(k)
\cap \F^{\rho}} F_{\rho}^{\s}\right|^2
\sum_{\rho_2}\left|b_Q^3\sum_{\s_2\in\tilde{\G}(\ell)\cap\G^{\rho_2}}
 G_{\rho_2}^{\s_2}\right|^2.$$
It follows by \bref{u1} that $A(k,\ell)$ is 
\begin{equation}\lesssim\min(|\tilde{\F}(k)|,|\tilde{\G}(\ell)|)\int b_Q^{12}\sum_{\stackrel{\rho,\s,\pi}{\s\in\F(k)
\cap\F^{
\rho}}}\sum_{\stackrel{\rho_2,\s_2,\pi_2}{\s_2\in\G(\ell)
\cap\G^{\rho_2}} }|F_{\rho}^{\s,\pi}|^2|G_{\rho_2}^{\s_2,\pi_2}|^2\phi_{\pi}^{-2}\phi_{\s}^{-2}
\phi_{\pi_2}^{-2}\phi_{\s_2}^{-2}.
\label{m1}\end{equation}

We claim next that for each pair $(\pi,\pi_2)$  we have
\begin{equation} \int b_Q^4\phi_{\pi}^2\phi_{\pi_2}^2\lesssim N^{\frac{d-2}{2}}.\label{pro1}\end{equation}
Namely, $\pi$ and $\pi_2$ each have one \lq\lq short" direction in which the width
 is $1$, and these directions lie in
$\Gamma_1$ and $\Gamma_2$ respectively, and are therefore transverse. It follows 
that $\pi\cap\pi_2$ is contained within a 
bounded distance of a $(d-2)$-plane, hence that \begin{equation}\int_{\pi\cap \pi_2}b_Q^4
\lesssim N^{\frac{d-2}{2}}.\label{pro2}\end{equation}
Estimate \bref{pro1}
is just a version of \bref{pro2} incorporating Schwartz tails, and is proved by 
estimating $\phi_{\pi}^2$ by an 
appropriate sum of constants times characteristic functions of translates of $\pi$ 
(and similarly with $\phi_{\pi_2}^2$) and then 
applying \bref{pro2} to the terms in the resulting series.

We now consider the terms in the sum \bref{m1}. For each pair $(\rho,\s,\pi)$ 
and $(\rho_2,\s_2,\pi_2)$ we have
\begin{eqnarray}&&
\hskip-.75in \int b_Q^{12}|F_{\rho}^{\s,\pi}|^2|G_{\rho_2}^{\s_2,\pi_2}|^2\phi_{\pi}^{-2}\phi_{\s}^{-2}
\phi_{\pi_2}^{-2}\phi_{\s_2}^{-2}\label{pro3}\\
&
\lesssim&\|b_Q^2
\phi_{\pi}^{-2}\phi_{\s}^{-1}F_{\rho}^{\s,\pi}\|_{\infty}^2\|b^2_Q\phi_{\pi_2}^{-2}\phi_{\s_2}^{-1}
G_{\rho_2}^{\s_2,\pi_2}
\|_{\infty}^2\int b_Q^4\phi_{\pi}^2\phi_{\pi_2}^2\nonumber\\
&\lesssim& N^{\frac{d-2}{2}}\|b_Q^2\phi_{\pi}^{-2}\phi_{\s}^{-1}F_{\rho}^{\s,\pi}\|_{\infty}^2\|b^2_Q
\phi_{\pi_2}^{-2}\phi_{\s_2}^{-1}
G_{\rho_2}^{\s_2,\pi_2}
\|_{\infty}^2 \nonumber\end{eqnarray}
by \bref{pro1}. It then follows that
\begin{eqnarray}
&&\label{pro11} \int b_Q^{12}|F_{\rho}^{\s,\pi}|^2|
G_{\rho_2}^{\s_2,\pi_2}|^2\\
&&\qquad\lesssim N^{\frac{d-2}{2}}\| \phi_{\s}^2b_Q^2\|^2_{\infty}\|
\phi_{\s_2}^2b_Q^2\|^2_{\infty}\|\phi_{\s}^{-3}\phi_{\pi}^{-2}F_{\rho}^{\s,
\pi}\|_{\infty}^2 \|\phi_{\s_2}^{-3}\phi_{\pi_2}^{-2}G_{\rho_2}^{\s_2,\pi_2}\|_{\infty}^2\nonumber
 \\
&&\qquad\lesssim 2^{-4k-4\ell}N^{\frac{d-2}{2}}\|\phi_{\s}^{-3}\phi_{\pi}^{-2}
F_{\rho}^{\s,
\pi}\|_{\infty}^2 \|\phi_{\s_2}^{-3}\phi_{\pi_2}^{-2}G_{\rho_2}^{\s_2,\pi_2}\|_{\infty}^2.\nonumber
\end{eqnarray}
The first inequality followed from  \bref{pro3} by rearranging some factors, 
and the second inequality followed from \bref{bref}. 

Using \bref{pro11} and Lemma 3.2 we may now bound  \bref{m1} by
$$\min(|\tilde{\F}(k)|,|\tilde{\G}(\ell)|)2^{-4k-4\ell}N^{-\frac{d+2}{2}}\sum_{\rho,\rho_2}\sum_{\stackrel{\s\in\tilde{\F}(k)
\cap\F^{\rho}}{\s_2\in
\tilde{\G}(\ell)\cap\G^{\rho_2}}}
\|\phi_{\s}^{-4}F_{\rho}^{\s}\|_2^2\|\phi_{\s_2}^{-4}
G_{\rho_2}^{\s_2}\|_2^2$$
which by definition of $h_1$ and $h_2$ is
$$\lesssim\min(|\tilde{\F}(k)|,|\tilde{\G}(\ell)|)2^{-4k-4\ell}N^{-\frac{d+2}{2}}
\cdot
 N^{d+1}|\tilde{\F}(k)|\,|\tilde{\G}(\ell)| h_1^2h_2^2.$$
We now use \bref{19.1}, and obtain a bound on \bref{m1} by
$$2^{-3k-3\ell}\mu\nu\min(2^k\mu,2^{\ell}\nu)N^{\frac{d}{2}}h_1^2h_2^2.$$
Summing over $k$ and $\ell$ gives the lemma.\enddemo
\enddemo

Fix $\e>0$ and then $\eta>0$ and partition $Q(N)$ in nonoverlapping\break $N^{1-\e}$-squares; 
the letter $R$ below will always denote one of these squares. We recall that $f$ and $g$ have $L^2$ norm $1$ and
are supported on $\Gamma_1^{(N)}$ and $\Gamma_2^{(N)}$ respectively.

\nonumproclaim{Lemma 3.5} On $Q(N)${\rm ,} for any $h_1$ and $h_2$ there are
decompositions
$$F_{h_1}=F_g+F_b \mbox{ and $F_b = \sum_RF_b^R$}$$
$$ G_{h_2}=G_g+G_b \mbox{ and $G_b = \sum_R\G_b^R$}$$
 where  $\su F_b^R\subset R${\rm ,} 
$\su G_b^R\subset R${\rm ,} and the following estimates hold{\rm .}
\begin{itemize}

\ritem{1.} $\int_{Q(N)}|F_gG_g|^2+|F_bG_g|^2+|F_gG_b|^2\lesssim N^{-\frac{d+2}{2}+C\e}$

\ritem{2.} For each $R$ we have $F_b^R=\alpha_R\widehat{f_R}$ and $G_b^R=\beta_R\widehat{g_R}${\rm ,}
 where $\alpha_R$ and $\beta_R$ are supported on $R$ and have $L^{\infty}$ norm $\leq 1${\rm ,}
 and   $f_R$ and $g_R$ are supported on the $N^{-(1-\e)}$\/{\rm -}\/neighborhoods of $\Gamma_1$
 and $\Gamma_2$ respectively{\rm ,} and  \end{itemize}
\begin{equation}\sum_R\|f_R\|_2^2+\|g_R\|_2^2\leq C_{\eta} N^{-\e+\eta}.\label{bref2}\end{equation}
\endproclaim

\demo{Proof} Let $\W=\F(h_1)$, $\B=\G(h_2)$. We can assume that both  $h_1$ and $h_2$ 
are greater than $N^{-B_1}$ where $B_1$ is a large dimension-dependent constant, since 
otherwise it is easy to check that the lemma is valid with $F_b$ and $G_b$ equal to zero.
 It follows that the cardinalities of $\W$ and $\B$ are bounded by $N^{B_2}$. 

We apply Lemma 1.1$'$  after rescaling by $N$; thus $\d$ in Lemma 1.1 is $N^{-\frac{1}{2}}$; 
and we also set $\e$ in Lemma 1.1$'$ equal to twice the present $\e$. 

For each $N^{1-\e}$-square
$R$ we then define
\begin{eqnarray*}
F_b^R&=&\left\{\begin{array}{ll}\sum_{\stackrel{\s\in\W}{\s\sim R}}F_{\rho}^{\s}&
\mbox{on $R$}\\0&\mbox{elsewhere}\end{array}\right.\\ \noalign{\vskip6pt}
 G_b^R&=&\left\{\begin{array}{ll}\sum_{\stackrel{\s\in\B}{\s\sim R}}G_{\rho}^{\s}&
\mbox{on $R$}\\0&\mbox{elsewhere.}\end{array}\right.\end{eqnarray*}
Define $F_b$ to be equal to $F_b^R$ on $R$ for each $R$ and similarly with $G_b$, 
and define $F_g=F-F_b$, $G_g=G-G_b$.

We will now show that
$$\int_{Q(N)}|F_gG|^2\lesssim N^{-\frac{d+2}{2}+C\e}.$$
Namely, fix a $\sqrt{N}$-square $Q$. Define $\mu$ to be the maximum on $Q$ of $\Phi_{\tilde{\F}}$, where $\tilde{\F}$
 is the tubes  $w\in\F(h_1)$ such that $w\not\sim Q$, 
and define $\nu$ to be the maximum on $Q$ of $\Phi_{\G(h_2)}$. By Lemma 3.4 we 
have 
$$\int_Q|F_gG|^2\lesssim h_1^2h_2^2\mu^2\nu N^{\frac{d}{2}}.$$
We now sum over $Q$ and use property~3  of Lemma~1.1. This gives
$$\int_{Q(N)}|F_gG|^2\lesssim N^{C\e}h_1^2h_2^2N^{\frac{d}{2}}|\F(h_1)|\,|\G(h_2)|$$
which is $\lesssim N^{C\e}\cdot N^{-\frac{d+2}{2}}$ by Lemma 3.3. We can clearly
 estimate $\int_{Q(N)}|F_gG|^2$ and $\int_{Q(N)}|F_gG_g|^2$ in the same way, and it 
follows that property~1  holds.

We have the following almost orthogonality estimate:
\begin{equation}\left\|\sum_{\s\sim R}F_{\rho}^{\s}\right\|_2^2\lesssim h_1^2N^{\frac{d+1}{2}}|
\{\s:\s\sim R\}|.\label{pro6}\end{equation}
Namely, for fixed $\rho$ we have
$$\Bigg\|\sum_{\stackrel{\s\sim R}{\s\in\F^{\rho}}}F_{\rho}^{\s}\Bigg\|_2^2\lesssim
\sum_{\stackrel{\s\sim R}{\s\in\F^{\rho}}}\|\phi_{\s}^{-4}F_{\rho}^{\s}\|_2^2
\lesssim h_1^2N^{\frac{d+1}{2}}|\{\s\in\F^{\rho}:\s\sim R\}|$$
where the first inequality follows from the Schwartz inequality since $\sum_{\s}\phi_{\s}^8\lesssim 1$ pointwise
 and the second follows from the definition 
of $h_1$. Lemma 3.1 implies that the functions $\sum_{\stackrel{\s\sim R}{\s\in\F^{\rho}}}F_{\rho}^{\s}$ are
 essentially orthogonal for different $\rho$
 and \bref{pro6} follows. 

Using Lemma 3.1 again we see that, on each fixed $N^{1-\e}$ square $R$, $F_b=
\sum_{\s\sim R}F_{\rho}^{\s}$ agrees with the Fourier transform of a function 
$f^0_R$ supported on the $\frac{C}{N}$-neighborhood of $\Gamma_1$.
We have
\begin{eqnarray*} \noalign{\vskip6pt}
\sum_R\|f^0_R\|_2^2=\sum_R\left\|\sum_{\s\sim
R}F_{\rho}^{\s}\right\|_2^2&\lesssim&h_1^2N^{\frac{d+1}{2}}\sum_R|\{\s:\s\sim R\}|\\
\noalign{\vskip6pt}
&\lesssim&h_1^2N^{\frac{d+1}{2}}
|\F(h_1)| N^{\eta}\\ \noalign{\vskip6pt}
&\lesssim& N^{\eta}. \\ \noalign{\vskip-6pt}
\end{eqnarray*}
The first inequality follows from \bref{pro6}, the second inequality follows
from property~1  of Lemma~1.1$'$ and the last inequality follows from Lemma~3.3. Now fix 
$R$ and take a suitable
Schwartz function $\kappa$ supported in $D(0,\frac{1}{2})$ and whose Fourier 
transform is $\geq 1$ on a large disc centered at the origin. Let 
$\kappa_R(x)=e^{ik\cdot x} N^{d(1-\e)}\kappa(N^{1-\e}x)$ for an appropriate
 $k$; if $k$ is chosen correctly then  $\widehat{\kappa_R}\geq 1$ on $R$. 
Define $\alpha_R=\frac{1}{\widehat{\kappa_R}}$ on $R$ and zero otherwise, 
and $f_R=\kappa_R\ast f_R^0$. Then $F_b^R=\alpha_R\widehat{f_R}$. To make the
 estimate \bref{bref2} we will use the following fact, which follows  from Schur's test:
\begin{quote}
If $s$ is a  function supported in  $D(0,r)$ with $\|s\|_{\infty}\leq|D(0,r)|^{-1}$ 
 and if supp$f$  intersects every disc of radius $r$ in measure $\leq\gamma r$, then
 $\|s\ast f\|_2\lesssim\gamma^{\frac{1}{2}}\|f\|_2$.
\end{quote}
 
We apply this with $s=\kappa_R$, $f=f_R^0$, $r\approx N^{1-\e}$, $\gamma\approx N^{-\e}$, 
which is justified since  $f_R^0$ is supported on the $\frac{C}{N}$-neighborhood of $\Gamma$.
 It follows that $\|f_R\|_2^2\lesssim N^{-\e}\|f_R^0\|_2^2$, so we have the part of \bref{bref2} 
which relates to $f$. We can of course treat $g$ the same way, so the proof is complete.\enddemo
 
We note also that the $L^2$ norms of $F_g, F_b, G_g$ and $G_b$ on $Q(N)$ are all bounded by a constant; it suffices
 to prove
this for $F_b$ and $G_b$, and for them it follows from \bref{bref2}.

\section{Proof of Theorem 1} 

\advance\eqcount by 21

We will use  a lemma from the previous work:
 
\proclaimtitle{\cite{B}, \cite{TV}}
\proclaim{Lemma}  In order to prove Theorem {\rm 1} it suffices to prove
 that
$$\int_{Q(N)}|\widehat{fd\s}\widehat{gd\s}|^p\lesssim N^{\gamma}\|f\|_2^p\|g\|_2^p$$
for fixed $p>1+\frac{2}{d}$ and $\gamma>0${\rm .}
\endproclaim

This   lemma originates in Section 4 of \cite{B}, 
and the  version stated above is a special case of Lemma 2.4 in part I of \cite{TV}. We also make
 a further reduction which follows  by the uncertainty principle in the usual way: it
 suffices to prove that if $f$ and $g$ are
functions with $L^2$ norm $1$ which are supported on the $\frac{1}{N}$-neighborhoods 
of  $\Gamma_1$ and $\Gamma_2$ respectively, then 
\begin{equation}\int_{Q(N)}|\hat{f}\hat{g}|^p\lesssim N^{-p+\gamma}
\label{411}\end{equation}
if $p>1+\frac{2}{d}$ and $\gamma>0$.

~

The rest of this section is the proof of  \bref{411}. 

Fix $p>1+\frac{2}{d}$  and let  $\phi(N)$  be the supremum of the quantity
\begin{equation}N^p \int_{Q(N)}|\hat{f}\hat{g}|^p\label{412}\end{equation}
over functions $f$ and $g$ with $L^2$ norm $1$ which are supported in the 
$\frac{1}{N}$-neighbor\-hoods of $\Gamma_1$ and $\Gamma_2$ respectively. Fix
 a sufficiently small $\e$ and then a much smaller $\eta$; we will  
show that 
\begin{equation}\phi(N)\leq C( 1+ N^{\eta}\phi(N^{1-\e}))\label{z1}
\end{equation}
for a suitable constant $C$.

Namely, choose $f$ and $g$ with $L^2$ norm $1$ so that the quantity \bref{412} 
is essentially maximized. Then choose $h_1$ and $h_2$ using the pigeonhole principle
 so that
$$\int_{Q(N)}|F_{h_1}G_{h_2}|^p\gtrsim(\log N)^{-2p}\phi(N)$$
where $F_{h_1}$ and $G_{h_2}$ were defined in Section 3. This is possible since it is easy
 to see that parameter values 
$h$ which are less than a high negative power of $N$ make a negligible contribution. Now
 apply Lemma 3.5 with this choice of $h_1$ and $h_2$. With notation as in Lemma 3.5 we
 have (by the triangle inequality)
$$\phi(N)\lesssim (\log N)^{2p}\int_{Q(N)}(|F_bG_g|^p+|F_gG_b|^p+|F_gG_g|^p)\; 
+\;(\log N)^{2p}\sum_R\int_R|\widehat{f_R}\widehat{g_R}|^p.$$
In the first term, we estimate the $L^p$ norm by the $L^1$ and $L^2$ norms using
 H\"older's inequality, and use that the $L^1$ norms of $F_gG_b$, $F_bG_g$ and $F_bG_b$ are bounded by a constant by
 the remark at the end of Section~3. In the
 second term, by definition of $\phi(N^{1-\e})$, we can estimate the integral over a
 fixed $R$ by 
$$N^{-(1-\e)p}\phi(N^{1-\e})\|f_R\|_2^p\|g_R\|_2^p.$$ Making these estimates we 
conclude that
\begin{eqnarray*}N^{-p}\phi(N)&\lesssim &(\log N)^{2p}\left(\int_{Q(N)}(|FG_g|^2+|F_gG|^2+|F_gG_g|^2)\right)^{p-1} 
\\&&+\ (\log N)^{2p} N^{-(1-\e)p}\phi(N^{1-\e})\sum_R\|f_R\|_2^p\|g_R\|_2^p.\end{eqnarray*}
We now use H\"older's inequality on the sum over $R$ and then insert the estimates in 
Lemma 3.5; this gives
\begin{eqnarray*}N^{-p}\phi(N)&\lesssim&(\log N)^{2p}\left(\int_{Q(N)}(|FG_g|^2+|F_gG|^2+|F_gG_g|^2)\right)^{p-1}
\;\\&& +\; 
(\log N)^{2p}N^{-(1-\e)p}\phi(N^{1-\e})\left(\sum_R\|f_R\|_2^{2}\right)^{\frac{p}{2}}
\left(\sum_R\|g_R\|^{2}_2\right)^{\frac{p}{2}}
\\&\lesssim&
(\log N)^{2p}N^{(p-1)(C\e-\frac{d+2}{2})}\\
&& + \; (\log N)^{2p}N^{-(1-\e)p}
\phi(N^{1-\e})\cdot N^{p(-\e+\eta) } .\end{eqnarray*}
The assumption $p>1+\frac{2}{d}$ implies that the exponent $p-\frac{d+2}{2}(p-1)$ 
is negative. We therefore obtain \bref{z1}, since we can replace $\eta$ by 
$\frac{\eta}{p+1}$, say. 

If $\gamma$ is given and if we take $\eta$ sufficiently small then estimate \bref{z1}
 implies by an obvious induction that $\phi(N)\lesssim N^{\gamma}$; thus we have
 proved \bref{411} and therefore Theorem 1.\enddemo

\vglue-20pt
\section{Further remarks}
\vglue-4pt
\advance\eqcount by 24

We will now prove the corollary which was stated in the introduction. We first rephrase
 it in a somewhat sharper form and in general dimensions. We will use mixed norms on 
$\Gamma$ splitting the $S^{d-2}$ and radial  variables:
$$\|f\|_{L^p(L^q)}\stackrel{def}{=}
\left(\int_{S^{d-2}}\left(\int_1^2|f(t\omega)|^qdt\right)^{\frac{p}{q}}d\omega\right)^{\frac{1}{p}}.$$

In the statement below, note that when $d=4$ the condition on $p$ reduces to $p>3$; 
by duality we obtain  a bound $\|\hat{f}\|_{L^p(L^2)}\lesssim\|f\|_p$ for any $p<\frac{3}{2}$, which clearly includes
  the result that was stated in the 
introduction.  When $d\geq 5$ the requirement that $p$ be larger than $2+\frac{4}{d}$ 
becomes significant so the statement becomes weaker.

\nonumproclaim{{C}orollary 1} Assume that  
$p>\max(2+\frac{4}{d}, 2+\frac{2}{d-2})${\rm .} Let $f$ be a function  on $\Gamma${\rm .}
 Then $\|\widehat{fd\s}\|_p\leq C_p\|f\|_{L^p(L^2)}${\rm .}
\endproclaim

\demo{Proof} This is the same as the proof of Theorem 2.2 in \cite{TVV}; see also \cite{TV}, where the
 rescaling maps for the cone employed below are used.

Fix a large number $N$ and a spherical cap $c\subset S^{d-2}$ centered at a point $e\in S^{d-2}$ with radius 
$N^{-1}$, i.e.\ $c=\{ \omega\in S^{d-2}:|\omega-e|<N^{-1}\}$. Let $\Gamma_c=\{x\in\Gamma:\frac{\ov{x}}{x_d}\in c\}$. 
Define $T_c$ to be the linear map such that $T(e,1)=(e,1)$, $T(e,-1)=N^2(e,-1)$ and $Ty=Ny$ if $y\in\R^d$ is orthogonal
 to $ (e,1)$ and $(e,-1)$. $T_c$ maps light rays to light rays and has the following metric properties:
\begin{equation} \mbox{det}T_N =N^d\label{deter}\end{equation}
and $T_c$ expands the distance between any two light rays contained in $c$ by a factor of roughly $N$, and
 roughly preserves distances on each individual such
light ray.

 Let $c_1$ and $c_2$ be two caps contained in $c$ separated
 by an amount comparable 
to $N^{-1}$ and let $f$ and $g$ be functions on $\Gamma$ with $L^{p}(L^2)$ norm 
$1$ which are supported on $\Gamma_{c_1}$ and $\Gamma_{c_2}$ respectively. Define $\tilde{f}d\s$ and $\tilde{g}d\s$
to be the measures obtained by pushing forward $fd\s$ and $gd\s$ by the map $T_c$. Then
$\tilde{f}$ and $\tilde{g}$ are functions on $\Gamma$ whose conical supports are at 
least a constant distance apart, and their $L^{p}(L^2)$ norms are comparable to 
$N^{-\frac{d-2}{p'}}$; hence their $L^2$ norms are at most $N^{-\frac{d-2}{p'}}$.
 Furthermore we 
have the formulae
$$\widehat{fd\s}=\widehat{\tilde{f}d\s}\circ T_c^{-1}$$
$$\widehat{gd\s}=\widehat{\tilde{g}d\s}\circ T_c^{-1}$$
and therefore, by \bref{deter} and Theorem 1, 
\begin{eqnarray*}\int|\widehat{fd\s}\widehat{gd\s}|^{\frac{p}{2}}&=&N^d\int|
\widehat{\tilde{f}
d\s}
\widehat{\tilde{g}d\s}|^{\frac{p}{2}}\\&\lesssim&N^{d-p\frac{d-2}{p'}}
=N^{d-(p-1)(d-2)}\end{eqnarray*}
for any $p>2+\frac{4}{d}$.
We now cover $S^{d-2}$ with caps $c_j$ of \lq\lq width" $N^{-1}$ as 
above 
and let $f_j$ be functions on $\Gamma$ with supp$f_j\subset\Gamma_{c_j}$. By 
applying the preceding estimate and summing
 over $j$ we obtain
$$\sum_{(j,k):{\rm dist}(c_j, c_k)\approx 
N^{-1}}\int|\widehat{\chi_jfd\s}\widehat{\chi_kfd\s}|^{\frac{p}{2}}
\lesssim N^{d-(p-1)(d-2)}\sum_j\|f_j\|_p^p.$$
The exponent of $N$ is negative  if $p>2+\frac{2}{d-2}$. The result now follows 
exactly as in \cite{TVV}, since the supports of the Fourier transforms of the 
functions $\widehat{\chi_jfd\s}\widehat{\chi_kfd\s}$ have finite overlap
if $N$ is fixed and  $\mbox{dist}(c_j, c_k)\approx 
N^{-1}$.\enddemo

We now consider the Mockenhaupt square function
$$SF(x)=\left(\sum_{\rho}|F_{\rho}|^2\right)^{\frac{1}{2}}$$
where $F=\widehat{fd\s}$ with $f$ supported on $\Gamma$, $f=\sum_{\rho} f_{\rho}$
 with $f_{\rho}$ supported in the sector $\rho$ of width about $N^{-\frac{1}{2}}$ 
and  $F_{\rho}=\widehat{f_{\rho}d\s}$. The following simple result appears natural in higher dimensions 
where the expected critical exponent is $2+\frac{2}{d-2}$; we do not consider the  
question of $L^4(\R^3)$ estimates  except to note that Theorem 1 can of course be substituted into the
 numerology  in \cite{TV}.

\nonumproclaim{{C}orollary 2} If $2\leq p\leq 2+\frac{4}{d}$ then there is an
 estimate
\begin{equation}\|F\|_p\lesssim N^{(\frac{1}{2}-\frac{1}{p})
\frac{d-2}{4}+\e}\|SF\|_p\label{lin}\end{equation}
for any $\e>0${\rm .}
\endproclaim

{\it Proof}.
We introduce a \lq\lq weaker" square function $\tilde{S}$ defined as follows: let
 $F=\widehat{fd\s}$ be as above, let $\Delta$ run through a covering of $\Gamma$ by discs of radius
 $N^{-\frac{1}{2}}$, suppose that $f_{\Delta}$ is 
supported in $\Delta$ and $f=\sum_{\Delta}f_{\Delta}, F_{\Delta}=\widehat{f_{\Delta}d\s}$ and
$$\tilde{S}F=\left(\sum_{\Delta}|F_{\Delta}|^2\right)^{\frac{1}{2}}$$

To prove \bref{lin} we consider first the \lq\lq bilinear" version;  in 
this version, one can prove a stronger
 result where $\tilde{S}F$ replaces $SF$. Thus we
let $f$ and $g$ as in Theorem 1 and $F=\widehat{fd\s}$, $G=\widehat{gd\s}$, and will show that
\begin{equation}\|FG\|_{\frac{p}{2}}\lesssim N^{(\frac{1}{2}-\frac{1}{p})
\frac{d-2}{2}+\e}\left(\|\tilde{S}F\|^2_p+\|\tilde{S}G\|^2_p\right).\label{bili1}\end{equation}

Namely,  we have
$$\|FG\|_{L^{\frac{p}{2}}(Q)}\lesssim N^{-\frac{1}{2}+\e}\|b_QF\|_{2}
\|b_QG\|_{2}$$
when $p>2+\frac{4}{d}$. This follows  by applying \bref{411} (with $N$
 replaced by $N^{\frac{1}{2}}$ and $Q(N)$ replaced by $Q$) to the functions 
$b_QF$ and $b_QG$. By interpolation with $L^2$ there is also an estimate
$$\|FG\|_{L^{\frac{p}{2}}(Q)}\lesssim N^{-\frac{d+2}{2}(\frac{1}{2}-\frac{1}{p})+\e}\|b_QF\|_{2}\|b_QG\|_{2}$$
when $2\leq p\leq 2+\frac{4}{d}$. The $b_QF_{\Delta}$'s are essentially orthogonal 
(their Fourier supports are essentially disjoint) so  we can 
estimate $\|b_QF\|_{2}$ by $\|b_Q\tilde{S}F\|_{2}$; using this and then H\"older's 
inequality we obtain
\begin{eqnarray*}\|FG\|_{L^{\frac{p}{2}}(Q)}&\lesssim& N^{-\frac{d+2}{2}
(\frac{1}{2}-\frac{1}{p})+\e}\|b_Q\tilde{S}F\|_{2}\|b_Q\tilde{S}G\|_{2}\\&\lesssim&N^{-\frac{d+2}{2}(\frac{1}{2}-\frac{1}{p})+\e}\cdot N^{\frac{d}{2}(1-\frac{2}{p}+\e)}\|b_Q\tilde{S}F\|_{p}\|b_Q\tilde{S}G\|_{p}\\
&\leq&N^{(\frac{1}{2}-\frac{1}{p})
\frac{d-2}{2}+\e}(\|b_Q\tilde{S}F\|^2_{p}+\|b_Q\tilde{S}G\|^2_{p}).\\
\noalign{\noindent 
Now take an $\ell^{\frac{p}{2}}$ sum over $Q$. Using the rapid decay of $b$ we
 obtain \bref{bili1}.
}\end{eqnarray*}

\phantom{again}
\vglue-30pt
The same argument clearly applies to $S$, so we also have
\begin{equation}\|FG\|_{\frac{p}{2}}\lesssim N^{(\frac{1}{2}-\frac{1}{p})
\frac{d-2}{2}+\e}(\|SF\|^2_p+\|SG\|^2_p).\label{bili}\end{equation}
In the case of $S$, since the  maps $T_c$ essentially take sectors contained in $\Gamma_c$ to sectors one can
 pass from the estimate \bref{bili} to the \lq\lq linear" one (i.e.\ \bref{lin}) by 
rescaling, just as in \cite{TV} or in the proof of Corollary 1.\hfill\qed
 
\demo{Further remarks} 1. It will be clear to the experts that one could also obtain a partial result on
  the 
(higher dimensional) cone multiplier/local smoothing problem using the estimate \bref{bili1} together with the usual technology
 as 
discussed for example in \cite{M} and an estimate for a  Nikodym type light ray maximal function, followed by 
another rescaling argument to pass from the bilinear to the linear estimate. 
We do not present this here because the estimate we have at present for the maximal function is rather  crude. 
 
\medbreak

2.  Let  $p_d=2+\frac{2}{d-2}$.
 It is natural to ask the following question: is there an
estimate
\begin{equation}\|\hat{f}\|_{L^{p'}(L^{p})}\lesssim\|f\|_{p'}\label{nn1}\end{equation}
provided $p>p_d$. One could also weaken this by asking instead for the estimate ($p\geq p_d$)
\begin{equation}\forall \e\exists C_{\e}:\|\hat{f}\|_{L^{p'}(L^{p})}\leq C_{\e}\l^{\e}\|f\|_{p'}\label{nn1.1}
\end{equation}
if supp$f\subset D(0,\l)$, $\l\geq 1$.

This statement would easily imply the restriction conjecture for the sphere
$S^{d-2}$. Namely, suppose that $f\in L^{p'}(\R^{d-1})$ with $p'$ as above and that
$f$ is supported in $Q(N)$, and apply \bref{nn1} to the function $f(\ov{x})e^{2\pi i x_d}\phi(\frac{x_d}{N})$ 
where $\phi$ is a
suitable bump function. (if one assumes instead \bref{nn1.1} then this argument still works using Tao's $\e$-removal 
lemma, see \cite{TV}  for example.)
Of course \bref{nn1} would also solve the  cone restriction problem, so it appears to be a 
natural common generalization. 

The statements \bref{nn1} or \bref{nn1.1} are also related to several other conjectures in the literature.
 For example, \bref{nn1.1} may be seen to be weaker than the \lq\lq Radon transform" conjecture 
in \cite{T}, and is therefore also weaker than the so-called local smoothing 
conjecture \cite{MSS}. We sketch the argument as follows: let $Rf$ be the Radon transform of $f$ restricted to
 the planes orthogonal to  light rays as discussed in \cite{T}; we will use the notation of that paper. Observe
 that the partial Fourier transform  of $Rf$ in the $s$ variable can be identified with the restriction of
$\hat{f}$ to the cone. Because of this, a rescaling argument followed by an application of the
Hausdorff-Young theorem in the $s$ variable shows, assuming  \cite[formula (33)]{T} (and that $p'\leq 2$!),  that 
if supp$f\subset D(0,\l)\subset \R^d$ then
$$\|\hat{f}\|_{L^{p'}(L^p)}\lesssim\l^{\frac{d-1}{p}+\alpha}\|f\|_{p'}.$$
Thus if \cite[(33)]{T} were true for all $\alpha>-\frac{d-1}{p}$ as is conjectured in \cite{T} then it would follow
 that  \bref{nn1.1} holds.

In the four dimensional case,  estimate \bref{nn1} is superficially similar 
to Corollary 1, the difference being  that the radial dependence is now $L^{p}$ instead of $L^2$, but since it would
 imply the restriction conjecture for
  $S^{2}$ it should not   be accessible using only \lq\lq soft" 
Kakeya information like our Lemma 1.1. 
\enddemo

\bigbreak\centerline{\bf Appendix: Estimates for the restricted X-ray transform} \bigbreak

The  motivation for this appendix was  to clarify the relationship
 between Lemma 1.1 and other
approaches that have been taken to the restriction of the X-ray transform to  the light rays -  see for example \cite{B}, 
\cite{GS}, \cite{GSW}, \cite{GU},   \cite{TV} and [18].  
This leads to a family of mixed norm  estimates which we formulate as Theorem A.1 below. 

Let ${\cal L}$ be the space of light rays with the integral defined by
$$\int_{{\cal L}}f(\ell)d\ell=\int_{S^{d-2}}\int_{Y(\omega)}f(\ell(y, \omega))dyd\omega.$$
Here $\ell(y,\omega)$ is the line through $y$ with direction $(\omega, 1)$, and
 $Y(\omega)$ is the hyperplane perpendicular to $\ell(0,\omega)$. We define mixed 
norms on $G$ by
$$\|f\|_{L^q(L^r)}=\left(\int_{S^{d-2}}\left(\int_{Y(\omega)}|f(\ell(y,
\omega))|^rdy\right)^{\frac{q}{r}}d\omega\right)^{\frac{1}{q}} . 
$$ We define the X-ray transform as an operator from
functions on
$\R^d$ to functions on ${\cal L}$ via
$$Xf(\ell)=\int_{\ell}f$$
and will be interested in estimates for $X$ from $L^p$ to $L^q(L^r)$.

We first discuss necessary conditions  in order to formulate a plausible conjecture; we omit details here.
Suppose that $X$ is bounded from $L^p$ to $L^q(L^r)$. Then dilations give the
 condition
\begin{equation}\frac{d}{p}-\frac{d-1}{r}=1.\label{7.1}\end{equation}
See e.g.\ \cite{C} and \cite{GSW}. Furthermore, the maps $T_c$ used in Section~5 give the condition
\begin{equation}\frac{d-2}{q}\geq\frac{d}{p}-\frac{d}{r}.\label{7.2}\end{equation}
Again see \cite{GSW}. Another condition can be obtained by considering the example
$f=\chi_E$ where $E$ is the $\d$-neighborhood of the cone segment $\Gamma$. This 
takes the form
\begin{equation}\frac{1}{p}\leq\frac{d}{2r}\label{7.3}\end{equation}

It is natural to expect that  \bref{7.1}, \bref{7.2}, \bref{7.3} are essentially also sufficient for boundedness. 
We will not consider endpoint questions and will therefore work locally. Index juggling  leads to the following

\vglue8pt 

\nonumproclaim{Plausible conjecture}  Let $p=q=\frac{d^2-2d+2}{d}$ and $r=\frac{d^2-2d+2}{2}${\rm .} 
Then
$X$ is bounded from the  Sobolev space $W^{p,\e}(Q(1))$ to $L^q(L^r)$ for any $\e>0${\rm .}
\endproclaim

\vglue8pt 

By $W^{p,\e}(Q(1))$ we mean functions supported in $Q(1)$ with
$$\|f\|_{p,\e}\stackrel{\rm def}{=}\|(1-\Delta)^{\frac{\e}{4}}f\|_p<\infty.$$

There is an obvious bound on $L^1$, namely, by Fubini's theorem
\begin{equation}\|Xf\|_{L^{\infty}(L^1)}=\|f\|_1.\label{4q5}\end{equation} 
Interpolating \bref{4q5} with 
the preceding conjecture we obtain the following conjectural bound on $L^p$.

\nonumproclaim{Plausible conjecture$_p$} $\phantom{{}_|}$\hskip-5pt  Assume that $ 1\leq p\leq\frac{d^2-2d+2}{d}${\rm
.}  Define $r$ via $\frac{d}{p}-\frac{d-1}{r}=1$ and $q$ via $\frac{d-2}{q}=\frac{d}{p}-\frac{d}{r}${\rm .} Then $X$ is 
bounded from  $W^{p,\e}(Q(1))$ to $L^q(L^r)$ for any $\e>0${\rm .} $\phantom{\sum^\int}$
\endproclaim
 
This would imply all local  $W^{p,\e}\rightarrow L^q(L^r)$ estimates with the given $p$ which are 
not ruled out by \bref{7.1} (in the local form where $\leq$ replaces $=$), \bref{7.2} and \bref{7.3}. 

We will prove the following:
 
\nonumproclaim{Theorem A.1} If $d=3$ or $d=4$ then the above conjectures are true{\rm .} If
 $d\geq 5$ then
the second conjecture is true on $L^p$ provided $p\leq\frac{d+1}{2}${\rm .}
\endproclaim

\demo{{R}emarks} 1. We note that $q$ and $r$ coincide when $p=\frac{d}{2}$, 
$q=r=d-1$, and that this case is covered
by our result. This is new except when $d=3$ (see below); it is 
analogous to the result of Drury \cite{D}  (see also \cite{OS} and \cite{C}) for the full X-ray transform. 
\medbreak

2. Consider the case $d=3$. In this case, the angular parameter $\omega$ runs
 over a one dimensional space and the restricted X-ray transform as defined here
 is a special case of the restricted X-ray transform associated to a \lq\lq rigid 
line complex"  \cite{GS}, \cite{GSW}. If $d=3$ and $q=r$, then the estimate in Theorem A.1 
is an estimate
from $W^{\frac{3}{2}, \e}$ to $L^2$. The latter estimate is known, 
actually in the sharper form where $\e=0$ -- cf.\ [18] (I thank Allan Greenleaf for this reference) and \cite{GS} -- and a
dual formulation of this same estimate is  used in \cite{TV}. However, Theorem A.1 
is new also in the three dimensional case if $p>\frac{3}{2}$. [Note added in proof: some higher-dimensional versions of these results have since been obtained in \cite{E}.]

\medbreak

3. It may be possible to obtain a scale invariant result (i.e.\ $\e=0$) by modifying the argument below, at 
least if one assumes strict inequality in \bref{7.2} and \bref{7.3} and ignores the three dimensional case, but we do not attempt that  here because the formulation of Lemma
 1.1  in the body of the paper is unsuitable for that purpose. We note though that our estimate on $W^{p\e}$ can
 immediately be \lq\lq upgraded" to a (local, of course) estimate on $L^p$ provided one assumes strict inequality 
in \bref{7.1}, \bref{7.2}, \bref{7.3}. This is because one can interpolate with the known fact  that $X$ is  bounded
 from a negative order $L^2$ Sobolev space to $L^2$.  We leave details to the  reader.

\medbreak

4.  A proof of the above conjectures  for the full range of $p$ in general dimensions
 has to be hard, since this would include a version of the Kakeya conjecture. Namely, if the first 
conjecture is true in $\R^d$, then a Kakeya set in $\R^{d-1}$ must have Minkowski
dimension at least $d+\frac{4}{d}-3$, as may be seen by applying the restricted $X$-ray
 bound to the indicator function of a cylinder over the $\d$-neighborhood of the Kakeya
 set. From this and  known arguments (namely the subadditivity of the minimal possible Minkowski dimension for 
a Kakeya set in $\R^n$ as a function of $n$) follows that
the first conjecture if true in all dimensions would imply that Kakeya sets have
 full Minkowski dimension. 
\enddemo

We will need the following numerical inequalities (trivial in principle, but we give proofs for
the reader's convenience). Here  $\theta\in[\frac{1}{2},1]$ (we emphasize that $\theta\geq \frac{1}{2}$) and the 
variables $x,y,a,b, a_j, b_k$ are nonnegative real numbers.
\begin{eqnarray}
&&\min(ax, by)^{\theta}\max(ax, by)^{1-\theta}\label{num1}\\
&&\qquad\qquad \leq\min(x, y)^{\theta}\max(x,
y)^{1-\theta}\max(a,b)^{\theta}\min(a,b)^{1-\theta}\nonumber\\ \noalign{\vskip6pt}
&&\min(\sum_j a_j,\sum_k b_k)^{\theta}\max(\sum_j a_j,\sum_k
b_k)^{1-\theta}\label{num2}\\
&&\qquad\qquad \leq\sum_{j,k}\min(a_j,b_k)^{\theta}\max( a_j, b_k)^{1-\theta}\nonumber\end{eqnarray}

\demo{Proofs} For \bref{num1} we may assume that $x\leq y$. If also $ax\leq by$, then
$$\min(ax, by)^{\theta}\max(ax, by)^{1-\theta}=a^{\theta}b^{1-\theta}\min(x, y)^{\theta}\max(x, y)^{1-\theta}$$
and \bref{num1} follows. If $ax\geq by$ then 
\begin{eqnarray*}
&&\hskip-1in\min(ax, by)^{\theta}\max(ax, by)^{1-\theta}\\
&=&(\frac{by}{ax})^{2\theta-1}a^{\theta}
x^{\theta}b^{1-\theta}y^{1-\theta}\\&\leq&a^{\theta}x^{\theta}b^{1-\theta}
y^{1-\theta}\\&=&\min(x, y)^{\theta}\max(x, y)^{1-\theta}\max(a,b)^{\theta}
\min(a,b)^{1-\theta}\end{eqnarray*}
since $a\geq b$.

For \bref{num2} we can assume $\sum_ja_j\leq\sum_kb_k=1$. In fact, we can assume in addition that $\sum_ja_j=1$.
 This follows from \bref{num1}: let $t=\sum_ja_j$ and consider the effect of replacing $a_j$ by $t^{-1}a_j$. The 
left side of \bref{num2} increases by a factor of $t^{-\theta}$, and
\bref{num1} implies the right side increases by at most this much.

The right side of \bref{num2} is smallest if $\theta=1$ so we are reduced to proving that
$\sum_ja_j=\sum_kb_k=1$ implies $\sum_{j,k}\min(a_j,b_k)\geq 1$. But 
$$\sum_j\sum_k\min(a_j,b_k)\geq\sum_j\min(a_j,\sum_kb_k)\geq\min(\sum_ja_j,\sum_kb_k)$$
so we are done.\enddemo

We start the proof of Theorem A.1 by giving  a convenient restatement of Lemma 1.1; this differs from
 Lemma 1.1 only in that the Schwartz tails have been discarded and entropy replaced by measure, and is therefore an immediate 
corollary of Lemma 1.1.

Let $\W$ and $\B$ be  $\d$-separated sets of white and black $\d$-tubes (thus they satisfy the 
transversality assumptions); 
assume each tube intersects the unit square. We let $\sim$ be the relation in Lemma~1.1 and will use the
 notation $w\sim x$ and $n_{\W}(Q)$ defined there. Let 
\begin{eqnarray*}
\Phi_{\W}(x)&=&\sum_{w\in\W}\chi_w(x),\;
\Phi_{\B}(x)=\sum_{b\in\B}\chi_b(x)\\
\noalign{\vskip6pt} 
\Phi^*_{\W}(x)&=&\sum_{\stackrel{w\in\W}{w\sim x}}\chi_w(x),\;
\Phi^*_{\B}(x)=\sum_{\stackrel{b\in\B}{b\sim x}}
\chi_b(x)\\
\noalign{\vskip6pt}
\tilde{\Phi}_{\W}&=&\Phi_{\W}-\Phi^*_{\W},\; \tilde{\Phi}_{\B}(x)=\Phi_{\B}-\Phi^*_{\B}.\end{eqnarray*}

\nonumproclaim{Lemma A.1}  The
following hold{\rm ,} where $C$ depends on $d$ only\/{\rm ;} the implicit constants also depend on 
$\e${\rm ,} and $Q$ runs over a partition of $Q(1)$ into $\d^{\e}$\/{\rm -}\/squares\/{\rm :}
\begin{itemize}
\ritem{1.}  $\sum_Q n_{\W}(Q)\lesssim |\W|\lo{5}${\rm .}

\ritem{2.}  $\sum_Q n_{\B}(Q)\lesssim |\B|\lo{5}${\rm .}

\ritem{3.} $|\{x\in Q(1):\tilde{\Phi}_{\W}(x)\geq\mu\mbox{ and }\Phi_{\B}(x)\geq\nu\}|\lesssim
\d^{-C\e} ({|\W|\,|\B|}/{\mu^2\nu})\d^d${\rm .}

\ritem{4.} $|\{x\in Q(1):\tilde{\Phi}_{\B}(x)\geq\nu\mbox{ and }\Phi_{\W}(x)\geq\mu\}|\lesssim
\d^{-C\e} ({|\W|\,|\B|}/{\mu\nu^2})\d^d$.\hfill\qed \end{itemize}
 
\endproclaim
 
The rough idea now is to regard properties~3  and 4  of Lemma~A.1 as a  \lq\lq virtual" 
$L^{\frac{3}{2}}$ to $L^3$  estimate
and to interpolate between this and an $L^1$ to $L^1$ estimate, namely the
 following:

\medbreak {\elevensc Lemma A.2.} $$|\{ x\in Q(1):\Phi_{\W}(x)\geq\mu \mbox{ and
}\Phi_{\B}(x)\geq\nu\}|\lesssim\d^{d-1}\min\left(\frac{|\W|}{\mu}, {|\B|}/{\nu}\right).$$

\demo{Proof} It is clear that $\|\sum_{w\in\W}\chi_w\|_{L^1(Q(1))}\lesssim
 |\W|\d^{d-1}$, hence the measure of the $\mu$-fold points is $\lesssim \frac{|\W|}
{\mu}\d^{d-1}$,  which implies the lemma.\enddemo

 Fix $\theta\in [\frac{1}{2}, 1]$ and define
\begin{eqnarray*}
\Psi_{\theta}&=&\min(\Phi_{\B},\Phi_{\W})^{\theta}\max(\Phi_{\B},\Phi_{\W})^{1-\theta}\\
S_{\theta}&=&\min(\Phi^*_{\B},\Phi^*_{\W})^{\theta}
\max(\Phi^*_{\B},\Phi^*_{\W})^{1-\theta}\\
T_{\theta}&=&(\tilde{\Phi}_{\B})^{\theta}\Phi_{\W}^{1-\theta}+(\tilde{\Phi}_{\W})^{\theta}\Phi_{\B}^{1-\theta}.
\end{eqnarray*}
We will use below that
\begin{equation}\Psi_{\theta}\lesssim S_{\theta}+T_{\theta}.\label{13.7}\end{equation}
This is a consequence of the numerical inequality 
\begin{eqnarray*}
\min(a+b, c+d)^{\theta}\max(a+b, c+d)^{1-\theta}&\lesssim&
a^{\theta}(c+d)^{1-\theta}+c^{\theta}(a+b)^{1-\theta}\\
&& +\ \min(b, d)^{\theta}\max(b, d)^{1-\theta}\end{eqnarray*}
which follows for example from \bref{num2}.  

We now estimate $T_{\theta}$ for appropriate $\theta$ by interpolation between Lemmas A.1 and A.2.

\nonumproclaim{Lemma A.3} Let $p$ and $q$ satisfy
$1\leq q\leq 3$ and  $\frac{1}{q}\geq\frac{2}{p}-1${\rm .} Let $\theta=\frac{1}{4}(3-\frac{1}{q})${\rm .}
 Then  
\begin{equation}\|\d^{d-2}T_{\theta}\|^q_{L^q(Q(1))}\lesssim \d^{-C\e}\left(\d^{2d-3}|\B|
\cdot\d^{2d-3}|\W|\right)^{\frac{q}{2p}}.\label{9.1}\end{equation}
\endproclaim

\demo{Proof} It suffices to consider the case where $\frac{1}{q}=\frac{2}{p}-1$ since the $\d$-separation 
implies that the quantity $(\d^{2d-3}|\B|
\cdot\d^{2d-3}|\W|)$ is  $\lesssim 1$.

Define $Y(\mu, \nu)$ to be the set where $ \tilde{\Phi}_{\W}\geq\mu$ and $\Phi_{\B}\geq\nu$.
Lemmas A.1 and A.2 give
$$|Y(\mu,\nu)|\lesssim\d^{-C\e}\min\left(\frac{|\W|\,|\B|}{\mu^2\nu}\d^d,\frac{|\W|}{\mu}\d^{d-1},\frac{|\B|}{\nu}
\d^{d-1}\right)$$
and  therefore also 
\begin{eqnarray*}
|Y(\mu,\nu)|&\lesssim&\d^{-C\e}\left(\frac{|\W|\,|\B|}{\mu^2\nu}\d^d\right)^{\frac{q}{p}
-1}\left(\frac{|\W|}{\mu}\d^{d-1}\right)^{1-\frac{q}{2p}}\left(\frac{|\B|}{\nu}\d^{d-1}\right)^{1
-\frac{q}{2p}}\\
&=&\d^{-C\e}\frac{(|\B|\,|\W|)^{\frac{q}{2p}}\d^{d-2+\frac{q}{p}}}
{\mu^{\theta q}\nu^{(1-\theta)q}}\end{eqnarray*}
where we used the value of $\theta$ to obtain the last line. Summing over dyadic levels for $\mu$ and $\nu$ between 
$1$ and a negative power of $\d$ 
gives
$$\|\Phi_{\B}^{1-\theta} (\tilde{\Phi}_{\W})^{\theta}\|_q^q\lesssim\d^{-C\e}(|\B|\,|\W|)^{\frac{q}{2p}}
\d^{d-2+\frac{q}{p}}.$$
This and the analogous estimate with the roles of $\B$ and $\W$ reversed imply
$$\|T_{\theta}\|_q^q\lesssim\d^{-C\e}(|\B|\,|\W|)^{\frac{q}{2p}}\d^{d-2+\frac{q}{p}}$$
which is equivalent to \bref{9.1} when $\frac{1}{q}=\frac{2}{p}-1$.
\enddemo
 
We will now pass to a similar estimate for  $\Psi_{\theta}$.  We will use a rescaling argument and
 induction on $\d$ like the final argument in \cite{Wr1} or \cite{Wr2}.
 The 
rescaling argument requires another relation between the exponents, which is essentially the dual
 relation to
\bref{7.1}. We remark at this point that the quantity which we need to estimate in order to prove Theorem A.1 is $\min(\Phi_{\B},\Phi_{\W})$ and not the slightly larger $\Psi_{\theta}$. It is possible that the 
 slightly stronger result obtained by considering $\Psi_{\theta}$ could prove useful, but the main reason we use $\Psi_{\theta}$ is
that the rescaling argument in the proof is difficult to carry out with $\min(\Phi_{\B},\Phi_{\W})$.

\nonumproclaim{Lemma A.4} Assume that $q\leq 3${\rm ,}
 $\frac{1}{q}\geq \frac{2}{p}-1${\rm ,} and $1\leq\frac{q}{p}\leq\frac{d}{d-1}${\rm .}
 Then for any $\e>0$ there is a constant $A_{\e}$  making
 the following estimate valid\/{\rm ;} here $\theta=\frac{1}{4}(3-\frac{1}{q})${\rm :}
\begin{equation}\|\d^{d-2}\Psi_{\theta}\|^q_{L^q(Q(1))}\leq A_{\e}\d^{-C\e}
\left(\d^{2d-3}|\B|\cdot\d^{2d-3}|\W|\right)^{\frac{q}{2p}}.\label{rad}\end{equation}
\endproclaim

{\it Proof}. We start with the following observation concerning rescaling.

\medbreak {\it Claim}. Suppose that $\d$ is small enough and that \bref{rad}
has been proved with $\d$ replaced by $\d^{1-\e}$. Let $Q$ be a $\d^{\e}$-cube, 
and let 
$\B$ and $\W$ be $\d$-separated sets of tubes. Then 
$$ \|\Psi_{\theta}\|_{L^q(Q)}^q\leq   \d^{\frac{C\e^2}{2}}\cdot A_{\e}\d^{-C\e}
(\d^{2d-3}|\W|\cdot\d^{2d-3}|\B|)^{\frac{q}{2p}}.$$
\medbreak

Namely, for each $w\in \W$ let $k(w)$ be the cardinality of the set of tubes $w_1\in\W$
 such that $w_1\cap Q$ is contained in the double of $w$; similarly for each $b\in\B$ 
let $k(b)$ be the cardinality of the set of tubes $b_1\in\B$ such that $b_1\cap Q$ 
is
 contained in the double of $b$. Notice that $k(w)$ and $k(b)$ are between $1$ and
 $\d^{-(d-2)\e}$. Let
$\W(\mu)=\{w\in\W:k(w)\in[\mu, 2\mu]\}$, $\B(\nu)=\{ b\in\B:k(b)\in [\nu, 2\nu]\}$, and (analogously to the earlier definitions) let
\begin{eqnarray*}\Phi_{\W}^{\mu}&=&\sum_{w\in W_{\mu}}\chi_w,\; \Phi_{\B}^{\nu}=\sum_{b\in B_{\nu}}\chi_b\\
 \Psi^{\mu\nu}_{\theta}&=&\min(\Phi^{\nu}_{\B},\Phi^{\mu}_{\W})^{\theta}\max(\Phi^{\nu}_{\B},
\Phi^{\mu}_{\W})^{1-\theta}.
\end{eqnarray*}
Then
\begin{equation}\Psi_{\theta}\leq\sum_{\mu,\nu}\Psi^{\mu\nu}_{\theta}
\label{13.11}
\end{equation} where the sum is over dyadic values of $\mu$ and $\nu$. This follows from \bref{num2}. 
\pagebreak
 
By \bref{13.11} and pigeonholing, there are values of $\mu$ and $\nu$ such that
 $$\|\Psi_{\theta}^{\mu\nu}\|_{L^q(Q)}\gtrsim(\log\frac{1}{\d})^{-2}\|\Psi_{\theta}\|_{L^q(Q)}.
$$ 
We assume without loss of generality that 
$\mu\geq \nu$.
Now   let $\ov{\B}$ (resp.\ $\ov{\W}$) be  subsets of $\B(\nu)$ (resp.\ $\W(\mu)$) 
which are maximal with respect to the following property:
\begin{itemize}
\item[$(\ast)$] If $b_1,b_2\in\ov{\B}$ (resp $\ov{\W}$), then $b_1\cap Q$ is not contained 
in the double of $b_2$. 
\end{itemize}

Let $\ov{\Phi}_{\B}$ (resp.\ $\ov{\Phi}_{\W}$) be the sums of the characteristic 
functions of the tubes of width
$C_0\d$ coaxial with the tubes in 
$\ov{\B}$ (resp.\ $\ov{\W}$), and $$\ov{\Psi}_{\theta}=\min(\ov{\Phi}_{\B},\ov{\Phi}_{\W})^{\theta}
\max(\ov{\Phi}_{\B},\ov{\Phi}_{\W})^{1-\theta}.$$
Then $\Phi_{\W}^{\mu}\lesssim\mu\ov{\Phi}_{\W}$ and $\Phi_{\B}^{\nu}\lesssim\nu\ov{\Phi}_{\B}$, pointwise on $Q$; 
this follows from maximality of $\ov{\W}$ and $\ov{\B}$ provided $C_0$ is large enough.  Hence also
$$\Psi_{\theta}^{\mu\nu}\lesssim\mu^{\theta}\nu^{1-\theta}\ov{\Psi}_{\theta}$$
by \bref{num1}. Taking $L^q$ norms we conclude that \begin{equation}\|\Psi_{\theta}^{\mu\nu}\|^q_{L^q(Q)}
\lesssim\mu^{\theta q}\nu^{(1-\theta)q}\|\ov{\Psi}_{\theta}\|^q_{L^q(Q)}.\label{6.1}\end{equation}
Furthermore property $(\ast)$ implies
\begin{equation} |\ov{\B}|\lesssim \nu^{-1}|\B(\nu)|,\; |\ov{\W}|\lesssim \mu^{-1}|\W(\mu)|.\label{6.2}
\end{equation}
We now dilate the situation by a factor $\d^{-\e}$. This maps $Q$ to a cube $Q'$ of side~$1$,
and maps $\ov{\B}$ and $\ov{\W}$ to $\d^{1-\e}$-separated families of $C_0\d^{1-\e}
$-tubes.
Accordingly we can apply the hypothesis that \bref{rad} holds at scale $\d^{1-\e}$. We 
conclude that
\begin{eqnarray*}
&&
\|\d^{(1-\e)(d-2)}\ov{\Psi}_{\theta}(\d^{\e}x)\|_{L^q(Q', dx)}^q\\
&&\qquad\qquad\lesssim\ A_{\e}\d^{-C\e(1-\e)}
(\d^{(2d-3)(1-\e)}|\ov{\W}|\cdot\d^{(2d-3)(1-\e)}|\ov{\B}|)^{\frac{q}{2p}}.\end{eqnarray*}  
Making the change of variables $x\rightarrow\d^{\e}x$ and factoring out the powers of 
$\d^{\e}$ we get
$$\d^{-d\e}\d^{-q(d-2)\e}\|\d^{d-2}\ov{\Psi}_{\theta}\|_{L^q(Q)}^q\lesssim A_{\e}\d^{-C\e(1-\e)}
\d^{-\frac{q}{p}(2d-3)\e}(\d^{2d-3}|\ov{\W}|\cdot\d^{2d-3}|\ov{\B}|)^{
\frac{q}{2p}}.$$
We now substitute in the estimates \bref{6.1} and \bref{6.2}, obtaining
\begin{eqnarray*}
&&\mu^{-\theta q}\nu^{-(1-\theta)q}\d^{-d\e}\d^{-q(d-2)\e}\|\d^{d-2}\Psi_{\theta}^{\mu\nu}\|_{L^q(Q)}^q\\
&&\qquad \lesssim
 (\mu\nu)^{-\frac{q}{2p}}
\cdot \d^{-\frac{q}{p}(2d-3)\e} \cdot
A_{\e}\d^{-C\e(1-\e)}(\d^{2d-3}|\W_{\mu}|\cdot\d^{2d-3}|\B_{\nu}|)^{\frac{q}{2p}}\end{eqnarray*} or equivalently 
\begin{eqnarray*}\|\d^{d-2}\Psi_{\theta}^{\mu\nu}\|_{L^q(Q)}^q&\lesssim 
&\mu^{q(\theta-\frac{1}{2p})}
\nu^{q(1-\theta-\frac{1}{2p})}
\cdot \d^{-\frac{q}{p}(2d-3)\e+d\e+q(d-2)\e}\\
&&\cdot A_{\e}\d^{-C\e(1-\e)}(\d^{2d-3}|\W_{\mu}|
\cdot\d^{2d-3}|\B_{\nu}|)^{\frac{q}{2p}} . \end{eqnarray*}
But $\nu\leq\mu\lesssim \d^{-(d-2)\e}$, and the exponents $q(\theta-\frac{1}{2p})$ and $q(1-\theta-\frac{1}{2p})$ are
 both nonnegative. Since $\frac{q}{p}\leq\frac{d}{d-1}$, a little juggling of indices shows that therefore 
$$ \mu^{q(\theta-\frac{1}{2p})}\nu^{q(1-\theta-\frac{1}{2p})}
\d^{-\frac{q}{p}(2d-3)\e+d\e+q(d-2)\e}\lesssim 1.$$ It follows 
that
$$\|\d^{d-2}\Psi_{\theta}^{\mu\nu}\|_{L^q(Q)}^q\lesssim \d^{C\e^2}\cdot A_{\e}
\d^{-C\e}(\d^{2d-3}|\W_{\mu}|\cdot\d^{2d-3}|\B_{\nu}|)^{\frac{q}{2p}}$$
and therefore
$$ \|\d^{d-2}\Psi_{\theta}\|_{L^q(Q)}^q\lesssim (\log\frac{1}{\d})^{2q}\d^{\frac{C\e^2}{2}}\cdot \d^{\frac{C\e^2}{2}}
A_{\e}
\d^{-C\e}(\d^{2d-3}|\W|\cdot\d^{2d-3}|\B|)^{\frac{q}{2p}}.$$
The  factor $(\log\frac{1}{\d})^{2q}\d^{\frac{C\e^2}{2}}$ 
is evidently small for small $\d$, so the proof of the claim is complete.

We assume now that \bref{rad} has been proved for parameter values $\d>\d_0$ for a 
certain $\d_0$ (the case where $\d$ is large is easy if $A_{\e}$ has been chosen
appropriately) and will prove it when $\d^{1-\e}>\d_0$. This will evidently establish
 the lemma. 

We use \bref{13.7}, and observe that a bound like \bref{rad} with $\Psi_{\theta}$ replaced by $T_{\theta}$ 
follows from Lemma A.3; the implicit constant in Lemma A.3 is small compared 
with $A_{\e}$ if $A_{\e}$ has been chosen appropriately. To estimate $S_{\theta}$,  
 subdivide $Q(1)$ in $\d^{\e}$-cubes $Q$. On each fixed $Q$ we can apply the claim 
to the restriction of $S_{\theta}$ to $Q$, replacing $\W$ by $\{w\in\W:w\sim Q\}$ and similarly with $\B$. 

We obtain for each $Q$
\begin{equation}\|\d^{d-2}S_{\theta}\|_{L^q(Q)}^q\leq \d^{\frac{C\e^2}{2}}
\cdot A_{\e}(\d^{2d-3}n_{\W}(Q)\cdot\d^{2d-3}n_{\B}(Q))^{\frac{q}{2p}}.\label{6.3}
\end{equation}
We now sum over $Q$ concluding that 
\begin{eqnarray*}\|\d^{d-2}S_{\theta}\|_{L^q(Q(1))}^q&\leq&\d^{\frac{C\e^2}{2}}\cdot A_{\e}\sum_Q(\d^{2d-3}n_{\W}(Q)
\cdot\d^{2d-3}n_{\B}(Q))^{\frac{q}{2p}}\\
&\leq&\d^{\frac{C\e^2}{2}}\cdot
 A_{\e}(\sum_Q\d^{2d-3}n_{\W}(Q))^{\frac{q}{2p}}
(\sum_Q\d^{2d-3}n_{\B}(Q))^{\frac{q}{2p}}\\
&\lesssim&\d^{\frac{C\e^2}{2}}
\cdot A_{\e}
(C\log\frac{1}{\d})^{\frac{5q}{p}}
(\d^{2d-3}|\W|)^{\frac{q}{2p}}
(\d^{2d-3}|\B|)^{\frac{q}{2p}}.\end{eqnarray*}
The three inequalities followed respectively from \bref{6.3}, from H\"older's inequality (recall that $q\geq p$)
and from
properties~1 and 2  of Lemma~A.1. The factor $\d^{\frac{C\e^2}{2}}\cdot(C\log\frac{1}{\d})^{5\frac{q}
{p}}$ is small for small $\d$; the result now follows by combining the last inequality 
with
the preceding bound for $\|T_{\theta}\|_{L^q(Q)}^q$.\enddemo

Lemma A.4 is our main estimate and the rest of the argument is basically just another rescaling argument. 
This is fairly routine, so we will omit some details. 
In order to carry  out the argument efficiently we first make some further definitions and remarks.

We define a map $X^*$ from functions on ${\cal L}$ to functions of $\R^d$ via 
$$X^*f(x)=\int_{S^{d-2}}f(\ell(x,\omega))d\omega.$$
This is easily seen to be the adjoint map to $X$. If $c$ is a spherical cap on $S^{d-2}$, then
define ${\cal L}_c$ to be the set of light rays $\ell\in {\cal L}$ whose direction is $(\omega, 1)$ for some
 $\omega\in c$. For given $p$ and $r$ and $\d$, and a set $Y\subset{\cal L}$, define
$${\cal E}_{\d}^{p,r}(Y)=\|\chi_{Y_{\d}}\|_{L^p(L^r)}$$
where $Y_{\d}$ is the $\d$-neighborhood of $Y$ (with respect to a smooth metric on ${\cal L}$).  

Next fix a cap $c$ centered at a point $e\in S^{d-2}$ with radius $\s$. The map $T_c$ in Section~5 takes light
 rays to light rays, so there is an action $T_{c}:{\cal L}\rightarrow{\cal L}$, which has the following metric
  properties:

\begin{itemize}

\item[(a)] If $Y\subset{\cal L}_c$ then $\|\chi_{T_cY}\|_{L^p(L^r)}\approx\s^{-\frac{d-2}{p}-\frac{d}{r}}\|\chi_{Y}\|_{p,r}$.

\item[(b)] If $Y\subset{\cal L}_c$ then ${\cal E}_{\d}^{p,r}(T_cY)\lesssim
\s^{-\frac{d-2}{p}-\frac{d}{r}}{\cal E}_{\s\d}^{p,r}(Y)$.
\end{itemize}

 Property (a) is proved as follows: within $c$, $T_c$ expands distances along $S^{d-2}$ by a factor $\s^{-1}$ (hence volumes
 by $\s^{-(d-2)}$), and if $\omega\in c$ then the action on the fiber $\{x\in\R^d:
x\bot(\omega, 1)\}$ expands volumes by roughly $$\overbrace{\s^{-1}\times\ldots\times\s^{-1}}^{d-2
\, times}\times \s^{-2},$$ i.e.\ by $\s^{-d}$. Thus $L^p(L^r)$ norms expand by $\s^{-\frac{d-2}{p}-\frac{d}{r}}$. 
Also property (b) follows from property (a) by observing that $T_{c}$ maps the $C\s\d$-neighborhood of $Y\subset{\cal
L}_c$ onto a
 set which includes the $\d$-neighborhood of 
$T_cY$.

Further  if $Y\subset {\cal L}_c$ then
\begin{equation} X^*\chi_Y(x)\approx \s^{d-2} X^*\chi_{T_cY}(T_cx).\label{13.1}\end{equation}
This follows from the definition of $X^*$ and the formula for volume expansion along $S^{d-2}$.

We will now rephrase Lemma A.4 using some of the preceding notation  and at the same time will replace it by a 
somewhat weaker result with a less cumbersome statement. 
 
\nonumproclaim{Lemma A.5} Let $Z\subset{\cal L}${\rm ,} let $C$ be a large constant and let $S$
 be a set of points in $\R^d$
 with the following properties\/{\rm :}
\begin{itemize}
\ritem{1.} The intersection of $S$ with the $\d$\/{\rm -}\/neighborhood of any given ray $\ell\in Z$ is contained in a cube of side
$1${\rm .}

\ritem{2.} If $x\in S$, then  there are two spherical caps $c_1$ and $c_2$ on $S^{d-2}$ with width $C^{-1}$ and whose
distance
 apart is at least $C^{-1}${\rm ,}  such that
$$\min(X^*(\chi_{{\cal L}_{c_1}\cap Z})(x), X^*(\chi_{{\cal L}_{c_2}\cap Z})(x))\geq\mu.$$
\end{itemize}

Then 
$$|S|\lesssim\d^{-\e}\mu^{-q}{\cal E}_{\d}^{p,r}(Z)^q$$
for any fixed $\e>0$, provided $q\leq 3$, $q\geq p\geq r$, and $\frac{1}{q}\geq\frac{2}{r}-1$, $\frac{q}{r}\leq
\frac{d}{d-1}$.
\endproclaim

\demo{Proof}  We first make a couple of reductions. First, it suffices to prove the lemma with   
assumption 1  replaced by the stronger assumption that $A\subset Q(1)$. This follows in a standard way using 
that $q\geq p\geq r$: if the result is proved  for $A$ contained in a  square of side $1$, then one can tile by 
such squares, take an $L^q$ sum over the squares and use   hypothesis~1. It then also suffices to prove Lemma~A.5 
when $p=r$, since ${\cal E}_{\d}^{p,r}(Z)$ increases with $p$ when $Z$ is contained in a fixed compact subset. In
 addition, it suffices by a simple covering argument to prove the lemma  assuming that the caps $c_1$ and $c_2$ in
hypothesis~2  
 are independent of $x$.

Now define $Z_i=Z\cap {\cal L}_{c_i}$, let $\W$ and $\B$ be  maximal $\d$-separated subsets of $Z_1$ and $Z_2$ 
respectively and (for each $w\in \W$)  let $D_w$ be the $\d$-disc in ${\cal L}$ centered at $w$. 
Then $X^*(\chi_{Z_1\cap D(w)})\lesssim \d^{d-2}\chi_w$, where on the right side
$\chi_w$ is the characteristic function of the $\d$-neighborhood of the line $w$. So $X^*(\chi_{Z_1})
\lesssim\sum_w\d^{d-2}\chi_{w}=\Phi_{\W}$, where $\Phi_{\W}$ is as in Lemma A.1. 
Accordingly  $\min(X^*(\chi_{{\cal L}_{c_1}\cap Z})(x), X^*(\chi_{{\cal L}_{c_2}\cap Z})(x))\lesssim\min(\Phi_{\B},\Phi_{\W})\leq\Psi_{\theta}$. The result now follows from Lemma A.4 using 
Tchebyshev's inequality and that $(\d^{2d-3}|\B|\d^{2d-3}|\W|)^{\frac{1}{2r}}\lesssim{\cal E}_{\d}^{r,r}(Z)$.
\enddemo
 
The point will now be that for appropriate values of the exponents the statement of Lemma A.5 is essentially 
invariant under the rescaling maps $T_c$.

\nonumproclaim{Lemma A.6} Assume that $q\leq 3${\rm ,} $q\geq p\geq r${\rm ,} $\frac{1}{q}\geq\frac{2}{r}-1${\rm ,}
 $\frac{q}{r}
\leq\frac{d}{d-1}${\rm ,} and
\begin{equation}\frac{d-2}{p}+\frac{d}{r}\leq d-2+\frac{d}{q}.\label{13.5}\end{equation}
Let $Y\subset{\cal L}${\rm .} Then
$$\|X^*\chi_Y\|_{L^q(Q(1))}\lesssim\d^{-\e}{\cal E}_{\d}^{p,r}(Y).$$
\endproclaim

{\it Proof}. A standard argument shows  that it will suffice to prove the corresponding 
distributional estimate \begin{equation}|\{x\in Q(1):X^*\chi_Y(x)\geq\l\}|\lesssim\d^{-\e}\l^{-q}{\cal E}_{\d}^{p,r}(Y)^q\label{13.6}\end{equation}
in the case where $\l$ is bounded below by a high power of $\d$, say
$$\l\geq\d^{\frac{B}{2}(d-2)}$$
where $B$ is a large constant depending on $d$. This is because of the $\d^{-\e}$ factors and the fact that very 
small values of $\l$ clearly make a negligible contribution.

 To prove \bref{13.6}, let $A=\{x\in Q(1):X^*\chi_Y(x)\geq\l\}$ and define $A_{\s}$ to be all points $x$ with 
the property that there are two $\s$-caps $c_1$ and $c_2$ on $S^{d-2}$ whose distance apart is between $\s$ and
 $C\s$ and such that  $X^*(\chi_{Y\cap {\cal L}_{c_i}}) \geq C^{-1}\d^{\e}\l$ for $i=1,2$. We claim that $\cup_{\s}A_{\s}\supset A$; the union 
is over dyadic $\s\geq\d^B$.

Namely, if $x\in A$, then take the smallest $\s$ such that $X^*(\chi_{Y\cap{\cal L}_{c}})(x)\geq(C\s)^{\frac{\e}{B}}\l$ 
for some cap $c$ of width $C\s$. (The lower bound on $\l$ implies that then $\s\geq\d^B$.) Consider a covering of $\s$
 by caps $c_i$ of width $\s$. The minimality of $\s$ implies that $X^*(\chi_{Y\cap{\cal L}_{c_i}})(x)$ is small 
compared with $X^*(\chi_{Y\cap{\cal L}_{c}})(x)$ for each fixed $i$. It follows that the contribution from a fixed finite number of the $c_i$'s is similarly small, and therefore there must be two  $c_i$'s, 
call them $c_1$ and $c_2$, which are at distance $\geq \s$ apart such that  $X^*(\chi_{Y\cap{\cal L}_{c_i}})(x)
\gtrsim \s^{\frac{\e}{B}}\l$ for $i=1$ and~$2$. This implies the claim.

 By pigeonholing we may now choose $\s$ so that 
\begin{equation}|A_{\s}|\geq\d^{\e}|A|.\label{4q4}\end{equation}
Cover $S^{d-2}$ with a family of $C\s$-caps $c_i$ with bounded overlap. This gives a further decomposition
$$A_{\s}=\cup_iA_{\s}^{c_i}$$
where $A_{\s}^{c_i}$ is the set of $x$ for which the two $\s$-caps $c_1$ and $c_2$ in the definition of $A_{\s}$ may
 be taken to be contained in $c_i$.

We now fix one of the $c_i$'s and apply Lemma A.5 to the sets $Z=\break T_{c_i}({\cal L}_{c_i}\cap Y)$ and
$S=T_{c_i}(A_{\s}^{c_i})$. Formula \bref{13.1} shows  that   hypothesis~2  is satisfied with $\mu\approx\s^{-(d-2)}\l$,
 and since $A\subset Q(1)$ and $T_{c}$ preserves lengths in the $(e,1)$ direction, one can easily see that hypothesis~1  is
also satisfied. It follows that
\begin{eqnarray*}|T_{c_i}(A^{c_i}_{\s})|&\lesssim&\d^{-\e}(\s^{-(d-2)}\l)^{-q}{\cal E}^{p,r}_{\d}(T_{c_i}({\cal L}_{c_i}
\cap Y))\\
&\lesssim&(\s^{-(d-2)}\l)^{-q}\s^{-q(\frac{d-2}{p}+\frac{d}{r})}{\cal E}^{p,r}_{\s\d}({\cal L}_{c_i}\cap Y)^q
\end{eqnarray*} 
by property (b) above. Thus, using also \bref{deter}
$$|A_{\s}|\lesssim \d^{-\e}\s^d(\s^{-(d-2)}\l)^{-q}\s^{-q(\frac{d-2}{p}+\frac{d}{r})}{\cal E}^{p,r}_{\s\d}
({\cal L}_{c_i}\cap Y)^q$$
which implies that
 $$|A_{\s}|\lesssim\d^{-\e} {\cal E}^{p,r}_{\s\d}({\cal L}_{c_i}\cap Y)^q$$
by the assumption \bref{13.5}.

Now observe that the $\s\d$-neighborhoods of the sets ${\cal L}_{c_i}\cap Y$ are essentially disjoint (no point
 $y\in {\cal L}$ belongs to more than a bounded number). Accordingly we can sum over $c_i$ to obtain 
$$|A_{\s}|\lesssim \d^{-\e} {\cal E}^{p,r}_{\s\d}( Y)^q.$$
We now use \bref{4q4} and the fact that ${\cal E}_{\e}^{p,r}$ increases with $\e$. The result follows.\hfill\qed

\demo{Proof of Theorem {\rm A.1}} Let $p,q,r$ be as in Theorem A.1. Because the statement is obtained by 
interpolation with \bref{4q5} we can assume that $p$ has its largest possible value, namely $\frac{5}{3}$ 
if $d=3$ and $\frac{d+1}{2}$ if $d\geq 4$. The
 following relations on the dual exponents will hold: $$\frac{p'}{r'}\leq\frac{d}{d-1},$$
$$\frac{d-2}{q'}+\frac{d}{r'}\leq d-2+\frac{d}{p'},$$
$$\frac{1}{p'}\geq\frac{2}{r'}-1,$$
 $$r'\leq q'\leq p'\leq 3.$$
Namely the first two are dual to \bref{7.1} and \bref{7.2} respectively. The third follows since $p\leq\frac{d+1}{2}$
 and $r$ is defined by \bref{7.1},  while the last is most easily checked by using the explicit values of $p$, $q$ 
and $r$. 
Thus Lemma A.6 is applicable and shows that
\begin{equation}\|X^*\chi_Y\|_{L^{p'}(Q(1))}\lesssim\d^{-\e}{\cal E}_{\d}^{q',r'}(Y).\label{13.15}\end{equation}

We now pass to the dual estimate.
If $f$ is supported in $Q(1)$ then we define  $X_{\d}f(\ell)=\d^{-(d-1)}\int_{\ell^{\d}}f$, where $\ell^{\d}$ is the 
tube of width $\d$  with axis $\ell$. 

Fix a nonnegative function $f$ supported in $Q(1)$ with $\|f\|_p=1$ and consider the quantity $\|X_{\d}f\|_{L^q(L^r)}$.
 By duality there is a function $g:{\cal L}\rightarrow\R$ such that $\|g\|_{L^{q'}(L^{r'})}=1$ and
$$\int_{\cal L} gX_{\d}f\gtrsim\|X_{\d}f\|_{L^q(L^r)}.$$
Since $X_{\d}f$ is roughly constant on $\d$-discs and since values of  $X_{\d}f$ which
are less than a high power of $\d$ make a negligible contribution to the norm,
we can then conclude that there is a function $g:{\cal L}\rightarrow\R$ with $\|g\|_{L^{q'}(L^{r'})}=1$,
 with
$$\int_{\cal L} gX_{\d}\chi_E\gtrsim\d^{\e}\|X_{\d}f\|_{L^q(L^r)}$$
and such that $g$ has the special form
\begin{equation} g=\mu\chi_Y\label{13.13}\end{equation}
where $\mu$ is a scalar, and the set $Y$ is a union of $\d$-discs. Note that this implies $\mu{\cal E}_{\d}^{q',r'}(Y)\lesssim 1$. We also let $\tilde{Y}$ be the corresponding union of $2\d$-discs.

Letting $g$ be as in \bref{13.13}, we have 
\begin{eqnarray*}\|X_{\d}f\|_{L^q(L^r)}&\lesssim&\d^{-\e}\int_{\cal L} \mu\chi_YX_{\d}f\\
&\lesssim&\d^{-\e}\int_{\cal L} \mu\chi_{\tilde{Y}}Xf\\
&=&\d^{-\e}\int_{\R^d}X^*(\mu\chi_{\tilde{Y}})f.\end{eqnarray*}
 Now apply \bref{13.15}
to $\chi_{\tilde{Y}}$  and  use H\"older's inequality, obtaining
\begin{equation}\|X_{\d}f\|_{L^q(L^r)}\lesssim\d^{-\e}\|f\|_p\label{10.2}\end{equation}
since  $\mu{\cal E}_{\d}^{q',r'}(\tilde{Y})\lesssim 1$. 

It remains to trade $\e$ derivatives for the $\d^{-\e}$ factor, which is done in the usual way. 
Suppose that $f$ has $W^{p\e}$-norm $1$ and has support in $Q(1)$. If $\phi$ is an appropriately chosen 
$C_0^{\infty}$ function and $\phi_j(x)=2^{dj}\phi(2^jx)$ then we can express  $f=g+\sum_j \phi_j\ast f_j$,
 where $\hat{g}$ has compact support,  
 and where $\sum_j2^{\eta j}\|f_j\|_p\lesssim\|f\|_{p\e}$ for small $\eta$. It follows using the smoothing 
effect of $\phi_j$ that
$$Xf\lesssim 1+\sum_jX_{2^{-j}}|f_j|$$
and now the theorem follows by applying \bref{10.2} with a small enough value of $\e$ to the terms in the series.\enddemo
\bye